 \documentclass[10pt]{article}
 \usepackage[top=1in,bottom=1in,left=1.5in,right=1.5in]{geometry}
  \usepackage{amsmath,amssymb}
  \usepackage{latexsym}
  \usepackage[dvips]{pstricks} 
   \usepackage{pst-node}
  \usepackage[pdftex]{graphicx}
  
  \newcommand{\C}{\mathbb{C}}
  
  \newcommand{\N}{\mathbb{N}}
  
  \newcommand{\R}{\mathbb{R}}

  \renewcommand{\b}{\mathbf{b}}

  \renewcommand{\c}{\mathbf{c}}
  \newcommand{\e}{\mathbf{e}}
  
  \newcommand{\bF}{\mathbf{F}}
  
  \newcommand{\gl}{\mathbf{GL}}
  
  \newcommand{\n}{\mathbf{n}}

  \newcommand{\br}{\mathbf{r}}

  \newcommand{\U}{\mathbf{U}}
  \newcommand{\uu}{\mathbf{u}}
  \newcommand{\vv}{\mathbf{v}}
  \newcommand{\V}{\mathbf{V}}
  \newcommand{\w}{\mathbf{w}}
  \newcommand{\W}{\mathbf{W}}
  \newcommand{\x}{\mathbf{x}}
  \newcommand{\X}{\mathbf{X}}
  \newcommand{\y}{\mathbf{y}}
  
  \newcommand{\z}{\mathbf{z}}
  \newcommand{\0}{\mathbf{0}}
  \newcommand{\1}{\mathbf{1}}

  \newcommand{\bG}{\mathbf{G}}

  \newcommand{\cB}{\mathcal{B}}
  \newcommand{\cC}{\mathcal{C}}

  \newcommand{\cF}{\mathcal{F}}
  \newcommand{\cG}{\mathcal{G}}
  \newcommand{\cH}{\mathcal{H}}
  
  \newcommand{\cJ}{\mathcal{J}}
  
  \newcommand{\cO}{\mathcal{O}}
  
  \newcommand{\cR}{\mathcal{R}}
  \newcommand{\cS}{\mathcal{S}}
  \newcommand{\cT}{\mathcal{T}}

  \newcommand{\cX}{\mathcal{X}}

  \newcommand{\rO}{\mathrm{O}}
  
  \newcommand{\uU}{\underline U}

  \newcommand{\lan}{\langle}
  \newcommand{\ran}{\rangle}
  \newcommand{\an}[1]{\lan#1\ran}
  
  \newcommand{\hs}{\hspace*{\parindent}}
  \newcommand{\proof}{\hs \textbf{Proof.\ }}
  
  \newcommand{\tr}{\mathop{\mathrm{tr}}\nolimits}

  \newcommand{\Gr}{\mathop{\mathrm{Gr}}\nolimits}
  
  \newcommand{\trans}{^\top}
  \newcommand{\qed}{\hspace*{\fill} $\Box$\\}

  \newcommand{\rC}{\mathrm{C}}

  \newcommand{\rS}{\mathrm{S}}

  \newcommand{\rank}{\mathrm{rank\;}}

  \newtheorem{theo}{\bfseries \hs Theorem}

  \newtheorem{lemma}[theo]{\bfseries \hs Lemma}

  \newtheorem{assum}[theo]{\bfseries \hs Assumption}
  \numberwithin{equation}{section} 

 \renewcommand{\span}{\mathrm{span}}

\begin{document}

 \title{Low-rank approximation of tensors}

 \author{
 Shmuel Friedland\footnotemark[1] and Venu Tammali\\
 Department of Mathematics, Statistics and Computer Science,\\
 University of Illinois at Chicago, \\
Chicago, Illinois 60607-7045,
 USA,\\
E-mail: \texttt{friedlan@uic.edu, vtamma2@uic.edu}
 }
 \renewcommand{\thefootnote}{\fnsymbol{footnote}}
 \footnotetext[1]{Supported by the NSF grant DMS-1216393 
 }

 \renewcommand{\thefootnote}{\arabic{footnote}}
 \date{December 10, 2014}
 \maketitle
 \begin{abstract}
 In many applications such as data compression, imaging or
 genomic data analysis, it is important to approximate a given
 tensor by a tensor that is sparsely representable. For
 matrices, i.e. $2$-tensors, such a representation can be
 obtained via the singular value decomposition, which allows to
 compute best rank $k$-approximations. 
 For very big matrices a low rank approximation using SVD is not computationally feasible.
 In this case different approximations are available.  It seems that variants of the CUR-decomposition are 
 most suitable.

 For $d$-mode tensors $\cT\in \otimes_{i=1}^d \R^{n_i}$, with
 $d>2$, many generalizations of the singular value decomposition
 have been proposed to obtain low tensor rank decompositions.
 The most appropriate approximation seems to be best $(r_1,\ldots,r_d)$-approximation, 
 which maximizes the $\ell_2$ norm of the projection of $\cT$ on $\otimes_{i=1}^d \U_i$,
 where $\U_i$ is an $r_i$-dimensional subspace $\R^{n_i}$.
 One of the most common methods is the \emph{alternating maximization method} (AMM).
 It is obtained by maximizing on one subspace $\U_i$, while keeping all other fixed, and alternating
 the procedure repeatedly for $i=1,\ldots,d$.  Usually, AMM will converge to a local best approximation.
This approximation is a fixed point of a corresponding map on Grassmannians.
We suggest a Newton method for finding the corresponding fixed point.
We also discuss variants of CUR-approximation method for tensors.

The first part of the paper is a survey on low rank approximation of tensors.
The second new part of this paper is a new Newton method for best $(r_1,\ldots,r_d)$-approximation.
We compare numerically different approximation methods.  

 \end{abstract}

 \noindent {\bf 2000 Mathematics Subject Classification.} 14M15, 15A18,
 15A69, 65H10, 65K10.
 
 \noindent {\bf Key words} Tensor, best rank one approximation, best $(r_1,\ldots,r_d)$-approximation, sampling,
 alternating maximization method, singular value decomposition, Grassmann manifold, fixed point, Newton method.

 \renewcommand{\thefootnote}{\arabic{footnote}}

 \section{Introduction}\label{S:intro}
 Let $\R$ be the field of real numbers.  Denote by $\R^\n=\R^{n_1\times\ldots\times n_d}:=\otimes_{i=1}^d \R^{n_j}$, where $\n=(n_1,\ldots,n_d)$,
 the tensor products of $\R^{n_1},\ldots,\R^{n_d}$.
 $\cT=[t_{i_1,\ldots, i_d}]\in \R^\n$
 is called a \emph{$d$-mode tensor}.  Note that the
 number of coordinates of $\cT$ is $N=n_1\ldots n_d$.
 A tensor $\cT$ is called a \emph{sparsely representable tensor} if it can
 represented with a number of coordinates that is much smaller than $N$.

 Apart from sparse matrices, the best known example of a sparsely representable $2$-tensor
 is a low rank approximation of a matrix
 $A\in \R^{n_1\times n_2}$.  A rank $k$-approximation of $A$ is
 given by $A_{\textrm{appr}}:=\sum_{i=1}^k \uu_i\vv_i\trans$, which can be identified
 with $\sum_{i=1}^k \uu_i\otimes \vv_i$.  To store
 $A_{\textrm{appr}}$ we need only the $2k$ vectors
 $\uu_1,\ldots,\uu_k\in\R^{n_1},\;\vv_1,\ldots,\vv_k\in\R^{n_2}$.
 A best rank $k$-approximation of  $A\in \R^{n_1\times n_2}$ can be computed via
 the  \emph{singular value decomposition}, abbreviated here as SVD, \cite{GolV96}.
 Recall that if $A$ is a real symmetric matrix, then the best rank $k$-approximation
 must be symmetric, and is determined by the spectral decomposition of $A$.

 The computation of the SVD requires $\mathcal {O}(n_1n_2^2)$
 operations and at least $\mathcal{O} (n_1 n_2)$ storage, assuming that $n_2\le n_1$. Thus,
 if the dimensions $n_1$ and $n_2$ are very large, then the
 computation of the SVD is often infeasible. In this case other
 type of low rank approximations are considered, see e.g.
 \cite{AM01,DV06,DKM06,FriKNZ06,FMMN11,FKV04,GTZ97}.

 For $d$-tensors with $d>2$ the situation is rather
 unsatisfactory. It is a major theoretical and computational
 problem to formulate good generalizations of low rank
 approximation for tensors and to give efficient
 algorithms to compute these approximations, see e.g.
 \cite{deLdV00,LMV00,ES09,FMMN11,FMPS13,Kho,KB09,MMD06,Ose11,OsT09,ZhaG01}.

 We now discuss briefly the main ideas of the approximation methods for tensors discussed in this paper. 
 We need to introduce (mostly) standard notation for tensors.  Let $[n]:=\{1,\ldots,n\}$ for $n\in\N$. 
 For $\x_i:=(x_{1,i},\ldots,x_{n_i,i})\trans \in \R^{n_i},i\in[d]$,
 the tensor $\otimes_{i\in[d]}\x_i=\x_1\otimes\cdots\otimes \x_d=\cX=[x_{j_1,\ldots,j_d}]\in\R^\n$
 is called a decomposable tensor, or rank one tensor if $\x_i\ne\0$ for $i\in [d]$.  That is, $x_{j_1,\ldots,j_d}=x_{j_1,1}\cdots x_{j_d,d}$
 for $j_i\in[n_i], i\in[d]$.  Let $\an{\x_i,\y_i}_i:=\y_i\trans\x_i$ be the standard inner product on $\R^{n_i}$ for $i\in[d]$.
 Assume that $\cS=[s_{j_1,\ldots,j_d}]$ and $\cT=[t_{j_1,\ldots,j_d}]$ are two given tensors in $\R^\n$.
 Then $\an{\cS,\cT}:=\sum_{j_i\in[n_i],i\in[d]} s_{j_1,\ldots,j_d}t_{j_1,\ldots,j_d}$ is the standard inner product on $\R^\n$.
  Note that
 \begin{eqnarray*}
 &&\an{\otimes_{i\in[d]}\x_i,\otimes_{i\in[d]}\y_i}=\prod_{i\in[d]}\an{\x_i,\y_i}_i,\\
 &&\an{\cT,\otimes_{i\in[d]}\x_i}=\sum_{j_i\in[n_i],i\in[d]} t_{j_1,\ldots,j_d}x_{j_1,1}\cdots x_{j_d,d}.
 \end{eqnarray*}
 The norm $\|\cT\|:=\sqrt{\an{\cT,\cT}}$ is called the Hilbert-Schmidt norm. (For matrices, i.e. $d=2$, it is called the Frobenius norm.)

 Let $I=\{1\le i_1<\cdots <i_l\le d\}\subset [d]$.  Assume that $\cX=[x_{j_{i_1},\cdots,j_{i_l}}]\in \otimes_{k\in[l]} \R^{n_{i_k}}$.
 Then the contraction $\cT\times \cX$ on the set of indices $I$ is given by:
\[\cT\times \cX=\sum_{j_{i_k}\in [n_{i_k}], k\in[l]} t_{j_1,\ldots,j_d} x_{j_{i_1},\ldots,j_{i_l}}\in\otimes_{p\in[d]\setminus I} \R^{n_p}.\]

 Assume that $\U_i\subset \R^{n_i}$ is a subspace of dimension $r_i$ with an orthonormal basis
 $\uu_{1,i},\ldots,\uu_{r_i,i}$  for $i\in[d]$.  Let $\U:=\otimes_{i=1}^d \U_i\subset \R^\n$. 
Then $\otimes_{i=1}^d \uu_{j_i,i}$, where $j_i\in[n_i],i\in[d]$, is an orthonormal basis in $\U$.
 We are approximating $\cT\in\R^{n_1\times \cdots\times n_d}$ by a tensor
 \begin{equation}\label{Sapprox}
 \cS=\sum_{j_i\in[r_i],i\in[d]} s_{j_1,\ldots,j_d} \uu_{j_1,1}\otimes\cdots  \otimes\uu_{j_d,d} \in \R^\n
 \end{equation}
 The tensor $\cS'=[s_{j_1,\ldots,j_d}]\in \R^{r_1\times \cdots\times r_d}$ is the \emph{core tensor} corresponding to $\cS$ in 
 the terminology of \cite{Tuc66}.

 There are two major problems: The first one is how to choose the subspaces $\U_1,\ldots,\U_d$.
 The second one is the choice of the core tensor $\cS'$.  
 Suppose we already made the choice of $\U_1,\ldots,\U_d$.
 Then
 $\cS=P_{\U}(\cT)$ is the orthogonal projection of $\cT$ on $\U$:
 \begin{equation}\label{PUfor}
 P _{\otimes_{i\in[d]}\U_i}(\cT)=\sum_{j_i\in[r_i],i\in[d]} \an{\cT,\otimes_{i\in[d]} \uu_{j_i,i}} \otimes_{i\in[d]} \uu_{j_i,i}.
 \end{equation}
 If the dimensions of $n_1,\ldots,n_d$ are not too big, then this projection can be explicitly carried out.
 If the dimension $n_1,\ldots,n_d$ are too big to compute the above projection, then one needs to introduce other approximations.
 That is, one needs to compute the core tensor $\cS'$ appearing in \eqref{Sapprox} accordingly.
 The papers  \cite{AM01,DV06,DKM06,FriKNZ06,FMMN11,FKV04,GTZ97,Kho,MMD06,Ose11,OsT09} essentially choose $\cS'$
 in a particular way. 
 
 We now assume that the computation of $P_{\U}(\cT)$ is feasible.  Recall that
 \begin{equation}\label{normprojfor}
 \|P_{\otimes_{i\in[d]}\U_i}(\cT)\|^2=\sum_{j_i\in [r_i],i\in[d]}|\an{\cT,\otimes_{i=1}^d \uu_{j_i,i}}|^2.
 \end{equation}
 The best $\br$-approximation of $\cT$, where $\br= (r_1,\ldots,r_d)$, in Hilbert-Schmidt norm is the
 solution of the minimal problem:
 \begin{equation}\label{bestapproxprb}
 \min_{\U_i,\dim \U_i=r_i, i\in[d]}\min_{\cX\in \otimes_{i\in[d]}^d\U_i} \|\cT-\cX\|.
 \end{equation}
 This problem is equivalent to the following maximum
 \begin{equation}\label{maxnormprb}
 \max_{\U_i,\dim \U_i=r_i, i\in[d]} \|P_{\otimes_{i\in[d]}\U_i}(\cT)\|^2.
 \end{equation}
 
 The standard \emph{alternating maximization method}, denoted by AMM, for solving \eqref{maxnormprb} is to solve the maximum problem, where
 all but the subspace $\U_i$ is fixed.  Then this maximum problem is equivalent to finding an $r_i$-dimensional subspace
 of $\U_i$ containing the $r_i$ biggest eigenvalues of a corresponding nonnegative definite matrix $A_i(\U_1,\ldots,\U_{i-1},\U_{i+1},\ldots,\U_d)\in
 \rS_{n_i}$.  Alternating between $\U_1,\U_2,\ldots,\U_d$ we obtain a nondecreasing sequence of norms of projections which converges to $v$.
 Usually, $v$ is a critical value of $\|P_{\otimes_{i\in[d]}\U_i}(\cT)\|$.   See \cite{LMV00} for details.

 Assume that $r_i=1$ for $i\in [d]$.  Then $\dim\U_i=1$ for $i\in[d]$.  In this case the minimal problem \eqref{bestapproxprb} is called a best rank one 
 approximation of $\cT$.  For $d=2$ a best rank one approximation of a matrix $\cT=T\in\R^{n_1\times n_2}$ is accomplished by the first left and right 
singular vectors and the  corresponding maximal singular value $\sigma_1(T)$.  The complexity of this computation is $\cO(n_1n_2)$ \cite{GolV96}.
Recall that the maximum \eqref{maxnormprb} is equal to $\sigma_1(T)$, which is also called the spectral norm $\|T\|_2$.
For $d>2$ the maximum \eqref{maxnormprb}  is called the spectral norm of $\cT$, and denoted by $\|\cT\|_\sigma$.
The fundamental result of Hillar-Lim \cite{HL13} states that the computation of $\|\cT\|_\sigma$ is NP-hard in general.
Hence the computation of best $\br$-approximation is NP-hard in general.

Denote by $\Gr(r,\R^n)$ the variety of all $r$-dimensional subspaces in $\R^n$, which is called Grassmannian or Grassmann manifold.
Let
\[ \1_d:=(\underbrace{1,\ldots,1}_d), \quad \Gr(\br,\n):=\Gr(r_1,n_1)\times \cdots \times \Gr(r_d,n_d).\]
Usually, the AMM for best $\br$-approximation  of $\cT$ will converge to a fixed point of a corresponding 
map $\bF_\cT:\Gr(\br,\n)\to \Gr(\br,\n)$.  This observation enables us to give a new Newton method for finding a best 
$\br$-approximation to $\cT$.  For best rank one approximation the map $\bF_\cT$ and the corresponding Newton method
was stated in \cite{FMPS11}.

This paper consists of two parts.  The first part surveys a number of common methods for low rank approximation methods of matrices and tensors.
We did not cover all existing methods here.  We were concerned mainly with the methods that the first author and his collaborators were studying, and 
closely related methods.  The second part of this paper is a new contribution to Newton algorithms related to best $\br$-approximations.
These algorithms are different from the ones given in \cite{ZhaG01,ES09,SL10}.  Our Newton algorithms are based on finding the fixed points corresponding 
to the map induced by the AMM.  In general its known that for big size problem, where each $n_i$ is big for $i\in[d]$ and $d\ge 3$, Newton methods are 
not efficient.  The computation associate the matrix of derivatives (Jacobian) is too expensive in computation and time.
In this case AMM or MAMM (modified AMM) are much more cost effective.  This well known fact is demonstrated in our simulations. 

We now briefly summarize the contents of this paper.  In \S\ref{S:SVD} we review the well known facts of singular value decomposition (SVD) and its use for 
best rank $k$-approximation of matrices.   For large matrices approximation methods using SVD are computationally unfeasible. 
Section \ref{S:lowrankapmat} discusses a number of approximation methods of matrices which do not use SVD.  The common feature of these methods is
sampling of rows, or columns, or both to find a low rank approximation.  The basic observation in \S\ref{Sub:rowsamp} is that, with high probability,
a best $k$-rank approximation of a given matrix based on a subspace spanned by the sampled row is with in a relative $\epsilon$ error to the best 
$k$-rank approximation given by SVD.  We list a few methods that use this observation.  However, the complexity of finding a particular $k$-rank approximation
to an $m\times n$ matrix is still $\cO(kmn)$, as the complexity truncated SVD algorithms using Arnoldi or Lanczos methods \cite{GolV96,LSY98}.
In \S\ref{Sub:CURap} we recall the CUR-approximation introduced in \cite{GTZ97}.  The main idea of CUR-approximation is to choose $k$ columns and rows 
of $A$, viewed as matrices $C$ and $R$, and then to choose a square matrix $U$ of order $k$ in such a way that $CUR$ is an optimal approximation of $A$.
The matrix $U$ is chosen to be the inverse of the corresponding $k\times k$ submatrix $A'$ of $A$.  The quality of CUR-approximation can be determined
by the ratio of $|\det A'|$ to the maximum possible value of the absolute  value of all $k\times k$ minors of $A$.  In practice one searches for this maximum
using a number of random choices of such minors.  A modification of this search algorithm is given in \cite{FMMN11}.  
The complexity of storage of $C,R,U$ is $\cO(k\max(m,n))$.  The complexity of finding the value of each entry of $CUR$  is $\cO(k^2)$.
The complexity of computation of $CUR$ is $\cO(k^2mn)$.  In \S\ref{S:fastapten} we survey CUR-approximation of tensors given in \cite{FMMN11}.
In \S\ref{S:prelbrap} we discuss preliminary results on best $\br$-approximation of tensors.  In \S\ref{sub:maxprob} we show that the minimum problem 
\eqref{bestapproxprb} is equivalent to the maximum problem \eqref{maxnormprb}.  In \S\ref{sub:singvten} we discuss the notion of singular tuples and 
singular valuesva tensor introduced in \cite{Lim05}.  In \S\ref{sub:apprprb} we recall the well known solution of maximizing 
$\|P_{\otimes_{i\in[d]}\U_i}(\cT)\|^2$ with respect to one subspace, while keeping other subspaces fixed.  In \S\ref{S:AMM} we discuss AMM for best 
$\br$-approximation and its variations.  (In  \cite{LMV00, FMPS13} AMM is called \emph{alternating least squares}, abbreviated as ALS.)
In \S\ref{sub:gendp} we discuss the AMM on a product space.  We mention a \emph{modified alternating maximization method} and
and \emph{2-alternating maximization method},  abbreviated as MAMM and 2AMM respectively, introduced in \cite{FMPS13}.  The
MAMM method consists of choosing the one variable which gives the steepest ascend of AMM.  2AMM consists of maximization with respect to a pair of variables, 
while keeping all other variables fixed.  In \S\ref{sub:AMMbrapr} we discuss briefly AMM and MAMM for best $\br$-approximations for tensors.
In \S\ref{sub:companbrap} we give the complexity analysis of AMM for $d=3, r_1\approx r_2\approx r_3$ and $n_1\approx n_2\approx n_3$.
In \S\ref{S:fixpt} we state a working assumption of this paper that AMM converges to a fixed point of the induced map, which satisfies certain
smoothness assumptions.  Under these assumptions we can apply the Newton method, which can be stated in the standard form in $\R^L$.
Thus, we first do a number of AMM iterations and then switch to the Newton method.
In \S\ref{sub:newtboneappr} we give a simple application of these ideas to state a Newton method for best rank one approximation.
This Newton method was suggested in \cite{FMPS11}.  It is different from the Newton method in \cite{ZhaG01}.
The new contribution of this paper is the Newton method which is discussed in \S\ref{S:bestbrapprox} and \S\ref{S:forDF}.
The advantage of our Newton method is its simplicity, which avoids the notions and tools of Riemannian geometry as for example in \cite{ES09,SL10}.
In simulations that we ran, the Newton method in \cite{ES09} was $20\%$ faster than our Newton method for best $\br$-approximation of $3$-mode tensors.
However, the number of iterations of our Newton method was $40\%$ less than in \cite{ES09}. 
In the last section we give numerical results of our methods for best $\br$-approximation of tensors. 
In \S\ref{S:NumRes} we give numerical simulations of our different methods applied to a real computer tomography (CT) data set (the so-called MELANIX data set of OsiriX).
The summary of these results are given in \S\ref{S:Conclusion}.

\section{Singular Value Decomposition}\label{S:SVD}
Let $A\in\R^{m\times n}\setminus\{0\})$.  We now recall well known facts on the SVD of $A$ \cite{GolV96}.
See \cite{Ste93} for the early history of the SVD.  Assume that $r=\rank A$.  Then there exist $r$-orthonormal sets of vectors 
$\uu_1,\ldots,\uu_r\in \R^m, \vv_1,\ldots,\vv_r\in \R^n$ such that we have:
\begin{eqnarray}\notag
&&A\vv_i=\sigma_i(A)\uu_i, \quad \uu_i\trans A=\sigma_i(A)\vv_i\trans, \quad i\in [r], \quad \sigma_1(A)\ge \cdots\ge \sigma_r(A)>0,\\
&&A_k=\sum_{i\in[k]} \sigma_i(A)\uu_i\vv_i\trans, \quad k\in [r], \; A=A_r. \label{SVDid}
\end{eqnarray} 
The quantities $\uu_i$, $\vv_i$ and $\sigma_i(A)$ are called the left, right $i$-th singular vectors and $i$-th singular value of $A$ 
respectively, for $i\in[r]$.
Note that $\uu_k$ and $\vv_k$ are uniquely defined up to $\pm 1$ if and only if $\sigma_{k-1}(A)>\sigma_k(A)>\sigma_{k+1}(A)$.
Furthermore for $k\in [r-1]$ the matrix $A_k$ is uniquely defined if and only if $\sigma_k(A)>\sigma_{k+1}(A)$.
Denote by $\cR(m,n,k)\subset \R^{m\times n}$
the variety of all matrices of rank at most $k$.  Then $A_k$  is a best rank-$k$ approximation of $A$:
\[\min_{B\in\cR(m,n,k)} \|A-B\|=\|A-A_k\|.\]
Let $\U\in \Gr(p,\R^m), \V\in\Gr(q,\R^n)$.  We identify $\U\otimes \V$ with 
\begin{equation}\label{defUVT}
\U\V\trans:=\span\{\uu\vv\trans,\quad \uu\in\U,\vv\in \V\}\subset \R^{m\times n}.
\end{equation}
Then $P_{\U\otimes \V}(A)$ is identified with the projection of $A$ on $\U\V\trans$ with respect to the standard
inner product on $\R^{m\times n}$ given by $\an{X,Y}=\tr XY\trans$.
Observe that 
\[\textrm{Range }A=\U_r^\star, \; \R^m=\U_r^\star\oplus(\U_r^\star)^\perp,  \textrm{ Range }A\trans =\V_r^\star,\;
\R^n=\V_r^\star\oplus(\V_r^\star)^\perp.\]
 Hence
\[P_{\U\otimes \V}(A)=P_{(\U\cap\U_r^\star)\otimes(\V\cap\V_r^\star)}(A)\Rightarrow \rank P_{\U\otimes \V}(A)\le 
\min(\dim \U,\dim \V,r).\]
Thus
\begin{eqnarray}\notag
&&\max_{\U\in\Gr(p,\R^m),\V\in\Gr(q,\R^n)}\|P_{\U\otimes \V}(A)\|^2=\|P_{\U_l^\star\otimes\V_l^\star}(A)\|^2=\sum_{j\in[l]}\sigma_j(A)^2,\\
&&\min_{\U\in\Gr(p,\R^m),\V\in\Gr(q,\R^n)}\|A-P_{\U\otimes \V}(A)\|^2=\|A-P_{\U_l^\star\otimes\V_l^\star}(A)\|^2=\sum_{j\in[r]\setminus[l]}\sigma_j(A)^2,
\notag\\
&&l=\min(p,q,r). \label{maxprojmat}
\end{eqnarray}

To compute $\U_l^\star,\V_l^\star$ and $\sigma_1(A),\ldots,\sigma_l(A)$  of a large scale
matrix $A$ one can use Arnoldi or Lanczos methods
\cite{GolV96,LSY98}, which are implemented
in the partial singular value decomposition.  This requires a substantial number
of matrix-vector multiplications with the matrix $A$ and thus
a complexity of at least $\cO(lmn)$.
\section{Sampling in low rank approximation of matrices}\label{S:lowrankapmat}
Let $A=[a_{i,j}]_{i=j=1}^{m,n}\in\R^{m\times n}$ be given.
Assume that $\b_1,\ldots,\b_m\in\R^n, \c_1,\ldots,\c_n\in\R^m$ are the columns of $A\trans$ and $A$ respectively.
($\b_1\trans,\ldots\b_m\trans$ are the rows of $A$.)
Most of the known fast rank $k$-approximation are using sampling of rows or columns of $A$, or both.
\subsection{Low rank approximations using sampling of rows}\label{Sub:rowsamp}
Suppose that we sample a set  
\begin{equation}\label{defsubI}
I:=\{1\le i_1<\ldots<i_s\le m\}\subset [m], \quad |I|=s,
\end{equation}
of rows $\b_{i_1}\trans,\ldots,\b_{i_s}\trans$, where 
$s\ge k$.  Let  $\W(I):=\span(\b_{i_1},\ldots,\b_{i_s})$.
Then with high probability the projection of the first $i$-th right singular vectors $\vv_i$ on $\W(I)$ is very close to $\vv_i$ for $i\in[k]$,
provided that $s\gg k$.  In particular, \cite[Theorem 2]{DV06} claims:
\begin{theo}\label{DesVemtheo} (Deshpande-Vempala)  Any $A\in\R^{m\times n}$ contains a subset $I$ of $s=\frac{4k}{\epsilon}+2k \log(k +1)$ rows such that there
is a matrix $\tilde A_k$ of rank at most $k$ whose rows lie in $\W(I)$ and
\[\|A-\tilde A_k\|^2\le (1+\epsilon)\|A-A_k\|^2.\]
\end{theo}  

To find a rank-$k$ approximation of $A$, one projects each row of $A$ on $\W(I)$ to obtain the matrix $P_{\R^m\otimes \W(I)}(A)$.
Note that we can view   $P_{\R^m\otimes \W(I)}(A)$ as an $m\times s'$ matrix, where $s'=\dim\W(I)\le s$.  Then find a best rank $k$-approximation
to  $P_{\R^m\otimes \W(I)}(A)$, denoted as $P_{\R^m\otimes \W(I)}(A)_k$.  
Theorem \ref{DesVemtheo} and the results of \cite{FKV04} yield that 
\[\|A-P_{\R^m\otimes \W(I)}(A)_k\|^2\le (1+\epsilon)\|A-A_k\|^2 + \eta \|A-P_{\R^m\otimes \W(I)}(A)\|^2.\]
Here $\eta$ is proportional to $\frac{k}{s}$, and can be decreased with more rounds of sampling.
Note that the complexity of computing  $P_{\R^m\otimes \W(I)}(A)_k$ is $\cO(ks'm)$.  The key weakness of this method is that
to compute $P_{\R^m\otimes \W(I)}(A)$ one needs $\cO(s'mn)$ operations.  Indeed, after having computed an orthonormal  basis of $\W(I)$,
to compute the projection of each row of $A$ on $\W(I)$ one needs $s'n$ multiplications.

An approach for finding low rank approximations of $A$ using random sampling of rows or columns is given in Friedland-Kave-Niknejad-Zare
\cite{FriKNZ06}.  Start with a random choice of $I$ rows of $A$, where $|I|\ge k$ and $\dim \W(I)\ge k$.  Find $P_{\R^m\otimes \W(I)}(A)$ and
$B_1:=P_{\R^m\otimes \W(I)}(A)_k$ as above.  Let $\U_{1}\in\Gr(k,\R^m)$ be the subspace spanned by the first $k$ left singular vectors of $B_1$.
Find $B_2=P_{\U_1\otimes \R^n}(A)$.  Let $\V_2\in\Gr(k,\R^n)$ correspond to the first $k$ right singular vectors of $B_2$.
Continuing in this manner we obtain a sequence of rank $k$-approximations $B_1,B_2,\ldots$.  It is shown in \cite{FriKNZ06} that 
$\|A-B_1\|\ge \|A-B_2\|\ge \ldots$ and $\|B_1\|\le \|B_2\|\le \ldots$.  One stops the iterations when the relative improvement of the approximation
falls below the specified threshold.  Assume that $\sigma_k(A)>\sigma_{k+1}(A)$.  Since best rank-$k$ approximation is a unique local minimum for the function 
$\|A-B\|, B\in \cR(m,n,k)$ \cite{GolV96},
it follows that in general the sequence $B_j, j\in\N$ converges to $A_k$.    It is straightforward to show that this algorithm is the AMM for low rank approximations 
given in \S\ref{sub:AMMbrapr}.  Again, the complexity of this method is $\cO(kmn)$.

Other suggested methods as \cite{AM01,DKM06,RV07,Sar06} seem to have the same complexity $\cO(kmn)$,
since they project each row of $A$ on some $k$-dimensional subspace of $\R^n$.
\subsection{CUR-approximations}\label{Sub:CURap}
Let 
\begin{equation}\label{defsubJ}
J:=\{1\le j_1<\ldots<j_t\le n\}\subset [n], \quad |J|=t,
\end{equation}
and $I\subset [m]$ as in \eqref{defsubI} be given.  Denote by $A[I,J]:=[a_{i_p,j_q}]_{p=q=1}^{s,t}\in \R^{s\times t}$.
CUR-approximation is based on sampling simultaneously the set of $I$ rows and $J$ columns of $A$ and  the approximation
matrix to $A(I,J,U)$ given by 
\begin{equation}\label{defCUR}
A(I,J,U):=CUR, \quad C:=A[[m],J], \; R:=A[I,[n]], \; U\in \R^{t\times s}. 
\end{equation}
Once the sets $I$ and $J$ are chosen the approximation $A(I,J,U)$ depends on the choice of $U$.  Clearly the row and the column spaces of $A(I,J,U)$
are contained in the row and column spaces of $A[I,[n]]$ and $A[[m],J]$ respectively.  
Note that to store the approximation $A(I,J,U)$ we need to store the matrices $C$, $R$ and $U$.  The number of these entries is $tm+sn+st$.
So if $n,m$ are of order $10^5$ and $s,t$ are of order $10^2$  the storages of $C,R,U$ can be done in Random Access Memory (RAM), while  
the entries of $A$ are stored in external memory.
To compute an entry of $A(I,J,U)$, which is an approximation of the corresponding  entry of $A$, we need $st$ flops.

Let $\U$ and $\V$ be subspaces spanned by the columns of $A[[m],J]$ and $A[I,[n]]\trans$ respectively.
Then $A(I,J,U)\in \U\V\trans$, see \eqref{defUVT}.

Clearly, a best CUR approximation is chosen by the least squares
principle:
\begin{eqnarray}\notag
&&A(I,J,U^\star)=A([m],J)U^\star A(I,[n]),\\ 
&&U^\star=\arg\min\{\|A-A([m],J)UA(I,[n])\|,\; U\in\R^{|J|\times |I|}\}.\label{bestCURap}
\end{eqnarray}
The results in \cite{FriT07} show that the least squares solution of \eqref{bestCURap} is given by:
\begin{equation}\label{Ustarfor}
U^\star=A([m],J)^\dagger A A(I,[n])^\dagger.
\end{equation}
Here $F^\dagger$ denotes the Moore-Penrose pseudoinverse of a matrix $F$.
Note that $U^\star$ is unique if and only if 
\begin{equation}\label{uniqUsol}
\rank A[[m],J]=|J|, \quad \rank A[I,[n]]=|I|.
\end{equation}
The complexity of computation of $A([m],J)^\dagger$ and  $A(I,[n])^\dagger$ are $\cO(t^2m)$ and $\cO(s^2n)$ respectively.
Because of the multiplication formula for $U^\star$, the complexity of computation of $U^\star$ is $\cO(stmn)$.

One can significantly improve the computation of $U$, if one tries to best fit the entires of the submatrix $A[I',J']$ 
for given subsets $I'\subset [m], J'\subset [n]$.  That is, let 
\begin{eqnarray}\notag
&&U^\star(I',J'):=\arg\min\{\|A[I',J']-A(I',J)UA(I,J')\|,\; U\in\R^{|J|\times |I|}\}=\\
&&A[I',J]^\dagger A[I',J'] A^\dagger [I,J'].
\label{UstarforI'J'}
\end{eqnarray}
(The last equality follows from \eqref{Ustarfor}.)  The complexity of computation of $U^\star(I',J')$ is $\cO(st|I'||J'|)$.

Suppose finally, that $I'=I$ and $J'$.  Then \eqref{UstarforI'J'} and the properties of the Moore-Penrose inverse yield that
\begin{equation}\label{UstarforIJ}
U^\star(I,J)=A[I,J]^\dagger, \quad A(I,J,U^\star(I,J))= B(I,J):=A[[m],J]A[I,J]^\dagger A[I,[n]].
\end{equation}
In particular $B(I,J)[I,J]=A[I,J]$.  Hence
\begin{eqnarray}\notag
&&A[[m],J]=B(I,J)[[m],J],\; A[I,[n]]=B(I,J)[I,[n]] \textrm{ if } |I|=|J|=k \textrm{ and } \det A[I,J]\ne 0,\\
&&B(I,J)=A[[m],J]A[I,J]^{-1} A[I,[n]].\label{invertAIJ}
\end{eqnarray}
The original $CUR$ approximation of rank $k$ has the form $B(I,J)$ given by \eqref{invertAIJ} \cite{GTZ97}.

Assume that $\rank A\ge k$.  We want to choose an approximation $B(I,J)$ of the form \eqref{invertAIJ}
which gives a good approximation to $A$.  It is possible to give an upper estimate for the maximum of the 
absolute values of the entries of $A-B(I,J)$ in terms of $\sigma_{k+1}(A)$, provided that $\det A[I,J]$ is relatively
close to 
 \begin{equation}
 \label{4.13.defmuk}
 \mu_k:=\max_{I\subset [m],J\subset[n],|I|=|J|=} |\det A[I,J]|>0.
 \end{equation}
 Let 
 \begin{equation}\label{4.13.defentwnrm}
 \|F\|_{\infty,e}:=\max_{i\in[m],j\in[n]} |f_{i,j}|, \quad F=[f_{i,j}]\in\R^{m\times n}.
 \end{equation}
 
 The results of \cite{GTZ97, GorT01} yield:
 \begin{equation}\label{4.13.basin}
 \|A- B(I,J)\|_{\infty,e}\le
 \frac{(k+1)\mu_k}{\det A[I,J]}\sigma_{p+1}(A).
 \end{equation}
(See also \cite[Chapter 4, \S13]{Fri15}.)

 To find $\mu_k$ is probably an NP-hard problem in general \cite{Fri13a}.
 A standard way to find $\mu_k$ is either a random search or greedy search \cite{GOSTZ10,Fri13a}.
 In the special case when $A$ is a symmetric positive definite matrix one can give the exact conditions when
 the greedy search gives a relatively good result \cite{Fri13a}.

 In the paper by Friedland-Mehrmann-Miedlar-Nkengla \cite{FMMN11} a good approximation $B(I,J)$ of the form  
 \eqref{invertAIJ} is obtained by a random search on the maximum value of the product of the significant singular values of $A[I,J]$.
 The approximations found in this way are experimentally better than the approximations found by searching for $\mu_k$.  
 \section{Fast approximation of tensors}\label{S:fastapten}
 The fast approximation of tensors can be based on several decompositions of tensors such as: Tucker decomposition \cite{Tuc66};
 matricizations of tensors, as unfolding and applying SVD one time or several time recursively, (see below); higher order singular value decomposition (HOSVD)
 \cite{deLdV00}, Tensor-Train decompositions \cite{Ose09,Ose11};  hierarchical Tucker decomposition \cite{Gra10,Hac11}.
 A very recent survey \cite{GKT13} gives an overview on this dynamic field.  In this paper we will discuss only the CUR-approximation.
 \subsection{CUR-approximations of tensors}\label{sub:CURap}
 Let $\cT\in\R^{n_1\times \ldots n_d}$.  In this subsection we denote the entries of $\cT$ as $\cT(i_1,\ldots,i_d)$ for $i_j\in [n_j]$ and $j\in [d]$.
 CUR-approximation of tensors is based on  matricizations of tensors.  
 The unfolding of $\cT$  in the mode $l\in [d]$ consists of rearranging the entries of $\cT$ as a matrix $T_l(\cT)\in \R^{n_l\times N_l}$, 
 where $N_l=\frac{\prod_{i\in [d]}n_i}{n_l}$. 
 More general, let $K\cup L=[d]$ be a partition of $[d]$ into two disjoint nonempty sets.  Denote $N(K)=\prod_{i\in K} n_i, N(L)=\prod_{j\in L} n_j$.
 Then unfolding $\cT$ into the two modes $K$ and $L$ consists of rearranging the entires of $\cT$ as a matrix $T(K,L,\cT)\in \R^{N(K)\times N(L)}$.

 We now describe briefly the $CUR$-approximation of $3$ and $4$-tensors as described by Friedland-Mehrmann-Miedlar-Nkengla \cite{FMMN11}.
 (See \cite{MMD06} for another approach to 

 \noindent
 CUR-approximations for tensors.)
 We start with the case $d=3$.  Let $I_i$ be a nonempty subset of $[n_i]$ for $i\in [3]$.  Assume that the following conditions hold:
 \[|I_1|=k^2, \quad |I_2|=|I_3|=k, \quad J:=I_2\times I_3\subset [n_2]\times [n_3].\]
 We identify $ [n_2]\times [n_3]$ with $[n_2n_3]$ using a lexicographical order.  We now take the CUR-approximation of $T_1(\cT)$ as given in \eqref{invertAIJ}:
 \[B(I_1,J)=T_1(\cT)[[n_1],J]T_1(\cT)[I_1,J]^{-1}T_1(\cT)[I_1,[n_2n_3]].\]
We view $T_1(\cT)[[n_1],J]$ as an $n_1\times k^2$ matrix.  For each $\alpha_1\in I_1$ we view $T_1(\cT)[\{\alpha_1\}, [n_2n_3]]$ as an $n_2\times n_3$ matrix $Q(\alpha_1):=[\cT(\alpha_1,i_2,i_3)]_{i_2\in[n_2], i_3\in[n_3]}$.  Let $R(\alpha_1)$ be  the $CUR$-approximation of $Q(\alpha_1)$ based on the sets $I_2,I_3$:
\[R(\alpha_1):=Q(\alpha_1)[[n_2],I_3] Q(\alpha_1)[I_2,I_3]^{-1} Q(\alpha_1)[I_2,[n_3]].\]
Let $F:=T_1(\cT)[I_1,J]^{-1}\in \R^{k^2\times k^2}$.  We view the entries of this matrix indexed by the row $(\alpha_2,\alpha_3)\in I_2\times I_3$
and column $\alpha_1\in I_1$.  We write these entries as $\cF(\alpha_1,\alpha_2,\alpha_3), \alpha_j\in I_j, j\in[3]$, 
which represent a tensor $\cF\in \R^{I_1\times I_2 \times I_3}$.
The entries of $Q(\alpha_1)[I_2,I_3]^{-1}$ are indexed by the row $\alpha_3\in I_3$ and column $\alpha_2\in I_2$.
We write these entries as $\cG(\alpha_1,\alpha_2,\alpha_3), \alpha_2\in I_2,\alpha_3\in I_3$, which represent a tensor $\cG\in \R^{I_1\times I_2 \times I_3}$.
Then the approximation tensor $\cB=[\cB(j_1,j_2,j_3)]\in \R^{n_1\times n_2\times n_3}$ is given by:
\begin{eqnarray*}
&&\cB(i_1,i_2,i_3)=\sum_{\alpha_1\in I_1,\alpha_j,\beta_j\in I_j, j=2,3}\\ 
&&\cT(i_1,\alpha_2,\alpha_3)\cF(\alpha_1,\alpha_2,\alpha_3)
\cT(\alpha_1,j_2,\beta_3)\cG(\alpha_1,\beta_2,\beta_3)\cT(\alpha_1,\beta_2,j_3).
\end{eqnarray*}

We now discuss a CUR-approximation for $4$-tensors, i.e. $d=4$.  Let $\cT\in \R^{n_1\times n_2\times n_3\times n_4}$ and
$K=\{1,2\},L=\{3,4\}$.  The rows and columns of $X:=T(K,L,\cT)\in \R^{(n_1n_2)\times (n_3n_4)}$ are 
indexed by pairs $(i_1,i_2)$ and $(i_3,i_4)$ respectively.
Let 
\[I_j\subset [n_j], \;  |I_j|=k,\; j\in [4], \quad  J_1:=I_1\times I_2,\;J_2:=I_3\times I_4.\]
First consider the CUR-approximation $X[[n_1n_2],J_2] X[J_1,J_2]^{-1} X[J_1,[n_3n_4]]$ viewed as tensor $\cC\in\R^{n_1\times n_2\times n_3\times n_4}$.  
Denote by $\cH(\alpha_1,\alpha_2,\alpha_3,\alpha_4)$ the $((\alpha_3,\alpha_4), (\alpha_1,\alpha_2))$
entry of the matrix  $X[J_1,J_2]^{-1}$.  So $\cH\in \R^{I_1\times I_2\times I_3\times I_4}$. Then
\[\cC(i_1,i_2,i_3,i_4)=\sum_{\alpha_j\in I_j,j\in[4]} \cT(i_1,i_2,\alpha_3,\alpha_4)\cH(\alpha_1,\alpha_2,\alpha_3,\alpha_4)\cT(\alpha_1,\alpha_2,i_3,i_4).\]

For $\alpha_j\in I_j,j\in[4]$ view  vectors $X[[n_1n_2],(\alpha_3,\alpha_4)]$ and $X[(\alpha_1,\alpha_2),[n_3n_4]]$ as matrices 
$Y(\alpha_3,\alpha_4)\in \R^{n_1\times n_2}$ and 
$Z(\alpha_1,\alpha_2)\in \R^{n_3\times n_4}$  respectively.  
Next we find the CUR-approximations to these two matrices using the subsets $(I_1,I_2)$ and $(I_3,I_4)$ respectively:
\begin{eqnarray*}
&&Y(\alpha_3,\alpha_4)[[n_1],I_2]Y(\alpha_3,\alpha_4)[I_1,I_2]^{-1}Y(\alpha_3,\alpha_4)[I_1,[n_2]],\\
&&Z(\alpha_1,\alpha_2)[[n_3],I_4]Z(\alpha_1,\alpha_2)[I_3,I_4]^{-1}Z(\alpha_1,\alpha_2)[I_3,[n_4]].
\end{eqnarray*}
We denote the entries of $Y(\alpha_3,\alpha_4)[I_1,I_2]^{-1}$ and $Z(\alpha_1,\alpha_2)[I_3,I_4]^{-1}$ by
 $\cF(\alpha_1,\alpha_2,\alpha_3,\alpha_4), \alpha_1\in I_1,\alpha_2\in I_2$
and $\cG(\alpha_1,\alpha_2,\alpha_3,\alpha_4), \alpha_3\in I_3,\alpha_4\in I_4$ respectively.
Then the CUR-approximation tensor $\cB$ of $\cT$  is given by:
\begin{eqnarray*}
&&\cB(i_1,i_2,i_3,i_4)=\sum_{\alpha_j,\beta_j\in I_j, j\in [4]} \cT(i_1,\beta_2,\alpha_3,\alpha_4)\cF(\beta_1,\beta_2,\alpha_3,\alpha_4)
\cT(\beta_1,i_2,\alpha_3,\alpha_4)\\
&&\cH(\alpha_1,\alpha_2,\alpha_3,\alpha_4)\cT(\alpha_1,\alpha_2,i_3,\beta_4)\cG(\alpha_1,\alpha_2,\beta_3,\beta_4)\cT(\alpha_1,\alpha_2,\beta_3,i_4).
\end{eqnarray*}

We now discuss briefly the complexity of the storage and computing an entry of the CUR-approximation $\cB$.
Assume first that $d=3$.  Then we need to store $k^2$ columns of the matrices $T_1(\cT)$, $k^3$ columns of $T_2(\cT)$ and $T_3(\cT)$, and
$k^4$ entries of the tensors $\cF$ and $\cG$.  The total storage space is $k^2 n_1+k^3(n_2+n_3)+2k^4$. To compute each entry of $\cB$ we need
to perform $4k^6$ multiplications and $k^6$ additions.

Assume now that $d=4$.  Then we need to store $k^3$ columns of $T_l(\cT), l\in[4]$ and $k^4$ entries of $\cF,\cG,\cH$.
Total storage needed is $k^3(n_1+n_2+n_3+n_4+3k)$.   To compute each entry of $\cB$ we need
to perform $6k^8$ multiplications and $k^8$ additions.

 \section{Preliminary results on best $\br$-approximation}\label{S:prelbrap}

\subsection{The maximization problem}\label{sub:maxprob}

 We first show that the best approximation problem \eqref{bestapproxprb} is equivalent to the maximum problem \eqref{maxnormprb},
 see \cite{LMV00} and \cite[\S10.3]{Hac12}.  The Pythagoras theorem yields that
 \[\|\cT\|^2 = \|P_{\otimes_{i=1}^d \U_i}(\cT)\|^2 + \|P_{(\otimes_{i=1}^d \U_i)^\perp}(\cT)\|^2,\quad
 \|\cT-P_{\otimes_{i=1}^d \U_i}(\cT)\|^2= \|P_{(\otimes_{i=1}^d \U_i)^\perp}(\cT)\|^2.\]
 (Here $(\otimes_{i=1}^d \U_i)^\perp$ is the orthogonal complement of $\otimes_{i=1}^d \U_i$ in $\otimes_{i=1}^d \R^{n_i}$.)
 Hence
 \begin{equation}\label{bestapproxid}
 \min_{\U_i\in\Gr(r_i,\R^{n_i}), i\in[d]}\|\cT-P_{\otimes_{i=1}^d \U_i}(\cT)\|^2=\|\cT\|^2-
\max_{\U_i\in \Gr(r_i,\R^{n_i}), i\in [d]} \|P_{\otimes_{i=1}\U_i}(\cT)\|^2.
 \end{equation}
 This shows the equivalence of \eqref{bestapproxprb} and \eqref{maxnormprb}.
 \subsection{Singular values and singular tuples of tensors}\label{sub:singvten}
 Let $\rS(n)=\{\x\in\R^n, \;\|\x\|=1\}$.   Note that one dimensional subspace $\U\in\Gr(1,\R^n)$ is $\span(\uu)$, where $\uu\in\rS(n)$.
 Let $\rS(\n):=\rS(n_1)\times\cdots\times\rS(n_d)$.  Then best rank one approximation problem 
 for $\cT\in\R^{\n}$ is equivalent to finding
 \begin{equation}\label{axprobbr1ap} 
 \|\cT\|_{\sigma}:=\max_{(\x_1,\ldots,\x_d)\in\rS(\n)} \cT\times (\otimes_{i\in[d]} \x_i).
 \end{equation}
 Let $f_\cT:\R^{\n}\to \R$ is given by $f_{\cT}(\cX)=\an{\cX,\cT}$.   Denote by $\rS'(\n)\subset \R^\n$ all rank one tensors of
 the form $\otimes_{i\in[d]} \x_i$, where $(\x_1,\ldots,\x_n)\in\rS(\n)$.  Let $f_\cT(\x_1,\ldots,\x_d):=f_\cT(\otimes_{i\in[d]}\x_i)$.
 Then the critical points of $f_\cT|\rS'(\n)$ are given by the Lagrange multipliers
 formulas \cite{Lim05}:
 \begin{equation}\label{defsingvalstup}
 \cT\times (\otimes_{j\in[d]\setminus\{i\}} \uu_j)=\lambda \uu_i,\quad i\in [d],\quad (\uu_1,\ldots,\uu_d)\in\rS(\n).
 \end{equation}
 One calls $\lambda$ and $(\uu_1,\ldots,\uu_d)$ a singular value and singular tuple of $\cT$.  For $d=2$ these are the singular values and 
 singular vectors of $\cT$.  The number of \emph{complex} singular values of a generic $\cT$  is given in \cite{FO14}.   This number increases
 exponentially with $d$.  For example for $n_1=\cdots=n_d=2$ the number of distinct singular values is $d!$.  (The number of real singular values
 as given by \eqref{defsingvalstup} is bounded by the numbers given in \cite{FO14}.)

 Consider first the maximization problem of $f_\cT(\x_1,\ldots,\x_d)$ over $\rS(\n)$ where we vary $\x_i\in\rS(n_i)$ and keep the other variables fixed.
 This problem is equivalent to the maximization of the linear form $\x_i\trans (\cT\times (\otimes_{j\in [d]\setminus\{i\}}\x_j))$. 
 Note that if $\cT\times (\otimes_{j\in [d]\setminus\{i\}}\x_j)\ne \0$ then this maximum is achieved for 
 $\x_i=\frac{1}{\|\cT\times (\otimes_{j\in [d]\setminus\{i\}}\x_j)\|}\cT\times (\otimes_{j\in [d]\setminus\{i\}}\x_j)$.

 Consider second the maximization problem of $f_\cT(\x_1,\ldots,\x_d)$ over $\rS(\n)$ where we vary $(\x_i, \x_j)\in\rS(n_i)\times\rS(n_j), 1\le i < j\le d$ 
 and keep the other variables fixed.  This problem is equivalent to finding the first singular value and the corresponding right and left singular vectors
 of the matrix $\cT\times (\otimes_{k\in[d]\setminus\{i,j\}}\x_k)$.  This can be done by using  use Arnoldi or Lanczos methods
 \cite{GolV96,LSY98}.  The complexity of this method is $\cO(n_in_j)$, given
 the matrix $\cT\times (\otimes_{k\in[d]\setminus\{i,j\}}\x_k)$. 
 \subsection{A basic maximization problem for best $\br$-approximation}\label{sub:apprprb}
 
 Denote by $\rS_n\subset \R^{n\times n}$ the space of real symmetric matrices.
 For $A\in \rS_n$ denote by $\lambda_1(A)\ge \ldots\ge\lambda_n(A)$ the eigenvalues of $A$ arranged in a decreasing order
 and repeated according to their multiplicities. 
 Let $\rO(n,k)\subset \R^{n\times k}$ be the set of all $n\times k$ matrices $X$ with $k$ orthonormal columns, i.e. $X\trans X=I_k$, where
 $I_k$ is $k\times k$ identity matrix.  We view $X\in \R^{n\times k}$ as composed of $k$-columns $[\x_1\ldots\x_k]$.
 The column space of $X\in \rO(n,k)$ corresponds to a $k$-dimensional subspace $\U\subset\R^n$.
 Note that $\U\in\Gr(k,\R^n)$ is spanned by the orthonormal columns of a matrix $Y\in O(n,k)$ if and  only if $Y=XO$, for some
 $O\in \rO(k,k)$.
 
 For $A\in \rS_n$ one has the Ky-Fan maximal characterization \cite[Cor. 4.3.18]{HJ88}
 \begin{equation}\label{KyFanchar}
 \max_{[\x_1\ldots\x_k]\in \rO(n,k)}\sum_{i=1}^k \x_i\trans A\x_i=\sum_{i=1}^k \lambda_i(A).
 \end{equation}
 Equality holds if and only if the column space of $X=[\x_1\ldots\x_k]$ is a subspace spanned by $k$ eigenvectors corresponding to $k$-largest 
 eigenvalues of $A$.  
 
 We now reformulate the maximum problem \eqref{maxnormprb} in terms of orthonormal bases of $\U_i, i\in [d]$.
 Let $\uu_{1,i},\ldots,\uu_{n_i,i}$ be an orthonormal basis of $\U_i$ for $i\in [d]$.  Then $\otimes_{i=1}^d \uu_{j_i,i}, j_i\in [n_i],i\in[d]$
 is an orthonormal basis of $\otimes_{i=1}^d \U_i$.  Hence
 \[\|P_{\otimes_{i=1}^d\U_i}(\cT)\|^2=\sum_{j_i\in [n_i],i\in[d]}\an{\cT,\otimes_{i=1}^d \uu_{j_i,i}}^2.\]
 Hence \eqref{maxnormprb} is equivalent to
 \begin{eqnarray}\label{maxnormprb1}
 &&\max_{[\uu_{1,i}\ldots\uu_{r_i,i}]\in\rO(n_i,r_i),i\in[d]} \sum_{j_i\in [n_i], i\in[d]}\an{\cT,\otimes_{i=1}^d \uu_{j_i,i}}^2=\\
 &&\max_{\U_i\in\Gr(r_i,\R^{n_i}),i\in[d]} \|P_{\otimes_{i=1}^d \U_i}(\cT)\|^2. \notag
 \end{eqnarray}
 
 A simpler problem is to find 
 \begin{eqnarray}\label{maxnormprbbas}
 &&\max_{[\uu_{1,i}\ldots\uu_{r_i,i}]\in\rO(n_i,r_i)} \sum_{j_i\in [n_i], i\in[d]}\an{\cT,\otimes_{i=1}^d \uu_{j_i,i}}^2=\\
 &&\max_{\U_i\in\Gr(r_i,\R^{n_i})} \|P_{\otimes_{i=1}^d \U_i}(\cT)\|^2, \notag
 \end{eqnarray}
 for a fixed $i\in [d]$.
 Let 
 \begin{eqnarray}\notag
 &&\uU:=(\U_1,\ldots,\U_d)\in\Gr(\br,\n),\\
 &&\Gr_i(\br,\n):=\Gr(r_1,n_1)\times\ldots\times
 \Gr(r_{i-1},n_{i-1})\times\Gr(r_{i+1},n_{i+1})\times\ldots\times\Gr(r_d,n_d),\notag\\
 &&\uU_i:=(\U_1,\ldots,\U_{i-1},\U_{i+1},\ldots,\U_d)\in \Gr_i(\br,\n),\notag\\
 &&A_i(\uU_i):=\sum_{j_l\in [r_l],l\in [d]\setminus\{i\}} 
 (\cT\times \otimes_{k\in [d]\setminus \{i\}} \uu_{j_k,k})(\cT\times \otimes_{k\in [d]\setminus \{i\}} \uu_{j_k,k})\trans.
\label{defAUi}
\end{eqnarray}

The maximization problem \eqref{maxnormprb1} reduces to the maximum problem \eqref{KyFanchar} with $A=A_i(\uU_i)$.
Note that each $A_i(\uU_i)$ is a positive semi-definite matrix.  Hence $\sigma_j(A_i(\uU_i))=\lambda_j(A_i(\uU_i))$ for $j\in [n_i]$.
Thus the complexity to find the first $r_i$ eigenvectors of $A_i(\uU_i)$ is $\cO(r_in_i^2)$.
Denote by $\U_i(\uU_i)\in\Gr(r_i,\R^{n_i})$ a subspace spanned by the first $r_i$ eigenvectors of $A_i(\uU_i)$.
Note that this subspace is unique if and only if
\begin{equation}\label{uniqcondUiAi}
\lambda_{r_i}(A_i(\uU_i))>\lambda_{r_i+1}(A_i(\uU_i)).
\end{equation}
Finally, if $\br=\1_d$ then each $A_i(\uU_i)$ is a rank one matrix.  
Hence $\U_i(\uU_i)=\span(\cT\times \otimes_{k\in [d]\setminus \{i\}} \uu_{1,k})$.  For more details see \cite{FM08}.
\section{Alternating maximization methods for best $\br$-approximation}\label{S:AMM}
\subsection{General definition and properties}\label{sub:gendp}
Let $\Psi_i$ be a compact smooth manifold for $i\in[d]$.  Define 
\[\Psi:=\Psi_1\times \cdots \times \Psi_d, \; 
\hat \Psi_i=(\Psi_1\times \cdots \times\Psi_{i-1}\times \Psi_{i+1}\times \cdots\times\Psi_d) \textrm{ for }i\in [d]. \]
We denote by $\psi_i$, $\psi=(\psi_1,\ldots,\psi_d)$ and $\hat\psi_i=(\psi_1,\ldots,\psi_{i-1},\psi_{i+1},\ldots,\psi_d)$
the points in $\Psi_i$, $\Psi$ and $\hat \Psi_i$ respectively.  Identify $\psi$ with $(\psi_i,\hat\psi_i)$ for each $i\in[d]$.
Assume that $f:\Psi\to\R$ is a continuous function with continuous first and second partial 
derivatives.  (In our applications it may happen that $f$ has discontinuities in first and second partial derivatives.)   
We want to find the maximum value of $f$ and a corresponding maximum point $\psi^\star$:
\begin{equation}\label{maxprobfPsi}
\max_{\psi\in\Psi} f(\psi)=f(\psi^\star).
\end{equation}
Usually, this is a hard problem, where $f$ has many critical points and a number of these critical points are local maximum points.
In some cases, as best $\br$ approximation to a given tensor $\cT\in\R^\n$, we can solve the maximization problem
with respect to one variable $\psi_i$ for any fixed $\hat\psi_i$:
\begin{equation}\label{maxprobPsii}
\max_{\psi_i\in \Psi_i} f((\psi_i,\hat\psi_i))=f((\psi_i^\star(\hat\psi_i),\hat\psi_i)),
\end{equation}
for each $i\in[d]$.

Then the \emph{alternating maximization method}, abbreviated as AMM, is as follows.  Assume that we start with an initial point
$\psi^{(0)}=(\psi^{(0)}_1,\ldots,\psi^{(0)}_d)=(\psi_1^{(0)},\hat \psi_1^{(0,1)})$.   Then we consider the maximal problem 
\eqref{maxprobPsii} for $i=1$ and $\hat\psi_1:=\hat \psi_1^{(0,1)}$.  This maximum is achieved for $\psi_1^{(1)}:=\psi_1^\star(\hat\psi_1^{(0,1)})$.
Assume that the coordinates $\psi_1^{(1)},\ldots,\psi_{j}^{(1)}$ are already defined for $j\in [d-1]$.  Let $\hat \psi_{j+1}^{(0,j+1)}:=
(\psi_1^{(1)},\ldots,\psi_j^{(1)},\psi_{j+2}^{(0)},\ldots,\psi_d^{(0)})$.  
Then we consider the maximum problem 
\eqref{maxprobPsii} for $i=j+1$ and $\hat\psi_{j+1}:=\hat \psi_{j+1}j^{(0,j+1)}$.  This maximum is achieved for 
$\psi_{j+1}^{(1)}:=\psi_{j+1}^\star(\hat\psi_{j+1}^{(0,j+1)})$.  Executing these $d$ iterations we obtain $\psi^{(1)}:=(\psi^{(1)}_1,\ldots,\psi^{(1)}_d)$.
Note that we have a sequence of inequalities:
\[f(\psi^{(0)})\le f(\psi_1^{(1)},\hat\psi_1^{(0,1)})\le f(\psi_2^{(1)},\hat\psi_2^{(0,2)})\le\cdots \le f(\psi_d^{(1)},\hat\psi_d^{(0,d)})=f(\psi^{(1)}).\]
Replace $\psi^{(0)}$ with $\psi^{(1)}$ and continue these iterations to obtain a sequence $\psi^{(l)}=(\psi_1^{(l)},\ldots,\psi_d^{(l)})$ for
$l=0,\ldots,N$.  Clearly,
\begin{equation}\label{nondecpropamm}
f(\psi^{(l-1)})\le f(\psi^{(l)})  \textrm{ for }l\in\N \Rightarrow \lim_{l\to\infty} f(\psi^{(l)})=M.
\end{equation}
Usually, the sequence $\psi^{(l)}, l=0,\ldots,$ will converge to $1$-semi maximum point $\phi=(\phi_1,\ldots,\phi_d)\in\Psi$.  
That is, $f(\phi)=\max_{\psi_i\in\Psi} f((\psi_i,\hat\phi_i))$ for $i\in [d]$.   
Note that if $f$ is differentiable at $\phi$ then $\phi$ is a critical point of $f$.
Assume that $f$ is twice differentiable at $\phi$.
Then  $\phi$ does not have to be a local maximum point \cite[Appendix]{FMPS13}.

The \emph{modified alternating maximization method}, abbreviated as MAMM, is as follows.
 Assume that we start with an initial point
$\psi^{(0)}=(\psi^{(0)}_1,\ldots,\psi^{(0)}_d)$.  Let $\psi^{(0)}=(\psi_i^{(0)},\hat \psi_i^{(0)})$ for $i\in[d]$.  Compute $f_{i,0}=
\max_{\psi_i\in\Psi_i} f((\psi_i,\hat\psi_i^{(0)}))$ for $i\in [d]$.  Let  $j_1\in\arg\max_{i\in[d]} f_{i,0}$.
Then $\psi^{(1)}=(\psi_j^\star(\hat\psi_{j_1}^{(0)}),\hat\psi_{j_1}^{(0)})$ and $f_1=f_{1,{j_1}}=f(\psi^{(1)})$.  
Note that it takes $d$ iterations to compute $\psi^{(1)}$.
Now replace $\psi^{(0)}$ with $\psi^{(1)}$ and compute $f_{i,1}=
\max_{\psi_i\in\Psi_i} f((\psi_i,\hat\psi_i^{(1)}))$ for $i\in[d]\setminus\{j_1\}$.   Continue as above to find $\psi^{(l)}$ for $l=2,\ldots,N$.
Note that for $l\ge 2$ it takes $d-1$ iterations to determine $\psi^{(l)}$.   
Clearly, \eqref{nondecpropamm} holds.  It is shown in \cite{FMPS13} that the limit $\phi$ of each convergent subsequence of the points $\psi^{(j)}$ is  
$1$-semi maximum point of $f$.

In certain very special cases, as for best rank one approximation, we can solve the maximization problem
with respect to any pair of  variables $\psi_i,\psi_j$ for $1\le i<j\le d$,
where $d\ge 3$ and all other variables are fixed.  Let 
\begin{eqnarray*}
&&\hat \Psi_{i,j}:=\Psi_1\times \cdots \times\Psi_{i-1}\times \Psi_{i+1}\times\cdots\times \Psi_{j-1}\times \Psi_{j+1}\times\cdots\times \Psi_d,\\
&&\hat\psi_{i,j}=(\psi_1,\ldots,\psi_{i-1},\psi_{i+1},\ldots,\psi_{j-1},\psi_{j+1},\ldots,\psi_d)\in \hat \Psi_{i,j},\quad \psi_{i,j}=(\psi_i,\psi_j)\in \Psi_i\times \Psi_j.
\end{eqnarray*}
View $\psi=(\psi_1,\ldots,\psi_d)$ as $(\psi_{i,j},\hat\psi_{i,j})$ for each pair $1\le i<j\le d$.
Then 
\begin{equation}\label{maxprobPsiij}
\max_{\psi_{i,j}\in \Psi_i\times \Psi_j} f((\psi_{i,j},\hat\psi_{i,j}))=f((\psi_{i,j}^\star(\hat\psi_{i,j}),\hat\psi_{i,j})).
\end{equation}
A point $\psi$ is called \emph{2-semi maximum point} if the above maximum equals to $f(\psi)$ for each pair $1\le i<j\le d$.

The \emph{$2$-alternating maximization method}, abbreviated here as 2AMM, is as follows.
Assume that we start with an initial point
$\psi^{(0)}=(\psi^{(0)}_1,\ldots,\psi^{(0)}_d)$. 
Suppose first that $d=3$.  Then we consider the maximization problem \eqref{maxprobPsiij} for $i=2, j=3$ and $\hat\psi_{2,3}=\psi_1^{(0)}$.
Let $(\psi_2^{(0,1)},\psi_3^{(0,1)})=\psi_{2,3}^\star(\psi_1^{(0)})$.  Next let $i=1,j=3$ and $\psi_{1,3}=\psi_2^{(0,1)}$.
Then $(\psi_1^{(0,2)},\psi_3^{(0,2)})=\psi_{1,3}^\star(\psi_2^{(0,1)})$.  Next let $i=1,2$ and $\hat\psi_{1,2}=\psi_3^{(0,2)}$.
Then $\psi^{(1)}=(\hat\psi_{1,2}^\star(\psi_3^{(0,2)}),\psi_3^{(0,2)})$.  Continue these iterations to obtain $\psi^{(l)}$ for $l=2,\ldots$.
Again, \eqref{nondecpropamm} holds.  Usually the sequence  $\psi^{(l)},l\in\N$ will converge to a 2-semi maximum point $\phi$.
For $d\ge 4$ the 2AMM can be defined appropriately see \cite{FMPS13}.

A \emph{modified $2$-alternating maximization method}, abbreviated here as M2AMM, is as follows.
Start with an initial point $\psi^{(0)}=(\psi^{(0)}_1,\ldots,\psi^{(0)}_d)$ viewed as $(\psi_{i,j}^{(0)},\hat\psi_{i,j}^{(0)})$,  for each pair $1\le i<j\le d$.
Let $f_{i,j,0}:=\max_{\psi_{i,j}\in \Psi_i\times \Psi_j} f((\psi_{i,j},\hat\psi_{i,j}^{(0)}))$.  Assume that 

\noindent
$(i_1,j_1)\in\arg\max_{1\le i <j\le d} f_{i,j,0}$.  Then $\psi^{(1)}=(\psi_{i_1,j_1}^\star(\hat\psi_{i_1,j_1}^{(0)}),\hat\psi_{i_1,j_1}^{(0)})$.
Let $f_{i_1,j_1,1}:=f(\psi^{(1)})$.
Note that it takes $d \choose 2$ iterations to compute $\psi^{(1)}$.   Now replace $\psi^{(0)}$ with $\psi^{(1)}$ and compute $f_{i,j,1}=
\max_{\psi_{i,j}\in\Psi_i\times \Psi_j} f((\psi_{i,j},\hat\psi_{i,j}^{(1)}))$ for all pairs $1\le i<j\le d$ except the pair $(i_1,j_1)$.   
Continue as above to find $\psi^{(l)}$ for $l=2,\ldots,N$.
Note that for $l\ge 2$ it takes ${d \choose 2}-1$ iterations to determine $\psi^{(l)}$.   
Clearly, \eqref{nondecpropamm} holds.  It is shown in \cite{FMPS13} that the limit $\phi$ of each convergent subsequence of the points $\psi^{(j)}$ is  
$2$-semi maximum point of $f$.
\subsection{AMM for best $\br$-aproximations of tensors}\label{sub:AMMbrapr}
Let $\cT\in \R^\n$.  For best rank one approximation one searches for the maximum of the function $f_\cT=\cT\times (\otimes_{i\in [d]}\x_i)$ 
on $\rS(\n)$, as in \eqref{axprobbr1ap}.  
For best $\br$-approximation one searches for the maximum of  the function $f_\cT=\|P_{\otimes_{i\in[d]}\U_i}\|^2$ on $\Gr(\br,\n)$, as in \eqref{maxnormprbbas}.
A solution to the basic maximization problem with respect to one subspace $\U_i$ is given in \S\ref{sub:apprprb}.

The AMM for best $\br$-approximation were studied first by  de Lathauwer-Moor-Vandewalle \cite{LMV00}.  The AMM is called in  \cite{LMV00}
 \emph{alternating least squares}, abbreviated as ALS.  A crucial problem is the starting point of AMM.  A high order SVD, abbreviated as HOSVD, 
for $\cT$, see \cite{deLdV00}, gives  
a good starting point for AMM.   That is, let $T_l(\cT)\in\R^{n_l\times N_l}$ be the unfolded matrix of $\cT$ in the mode $l$, as in \S\ref{sub:CURap}. 
Then $\U_l$ is the subspace spanned by the first $l$-left singular vectors of $T_l(\cT)$.  The complexity of computing $\U_l$ is $\cO(r_lN)$, where
$N=\prod_{i\in[d]} n_i$.  Hence for large $N$ the complexity of computing partial HOSVD is high.  Another approach is to choose the starting subspaces
at random, and repeat the AMM for several choices of random starting points. 

MAMM for best rank one approximation of tensors was introduced by Friedland-Mehrmann-Pajarola-Suter in \cite{FMPS13} by the name
\emph{modified alternating least squares}, abbreviated as MALS.  2AMM for best rank one approximation was introduced in
\cite{FMPS13} by the name \emph{alternating SVD}, abbreviated as ASVD.
It follows from the observation that $A:=\cT\times (\otimes_{l\in[d]\setminus\{i,j\}}\x_l)$
is an $n_i\times n_j$ matrix.  Hence the maximum of the bilinear form $\x\trans A\y$ on $\rS((n_i,n_j))$ is $\sigma_1(A)$.  See \S\ref{S:SVD}. 
M2AMM was introduced in \cite{FMPS13} by the names MASVD.

We now introduce the following variant of 2AMM for best $\br$-rank approximation, called \emph{$2$-alternating maximization method variant} and 
abbreviated as 2AMMV.  Consider the maximization problem for a pair of variables as in \eqref{maxprobPsiij}.
Since for $\br\ne \1_d$ we do not have a closed solution to this problem, we apply the AMM for two variables $\psi_i$ and $\psi_j$, while keeping $\hat \psi_{i,j}$ 
fixed.  We then continue as in 2AMM method.
\subsection{Complexity analysis of AMM for best $\br$-approximation}\label{sub:companbrap}
Let $\uU=(\U_1,\ldots,\U_d)$.  Assume that 
\begin{eqnarray}\notag
&&\U_{i}=\span(\uu_{1,i},\ldots,\uu_{r_i,n_i}), \quad \U_i^\perp=\span(\uu_{r_i+1},\ldots, \uu_{n_i,i}),\\
&&\uu_{j,i}\trans \uu_{k,i}=\delta_{j,k},\;j,k\in [n_i], \quad i\in [d].\label{basisUi}
\end{eqnarray}
For each $i\in[d]$ compute the symmetric positive semi-definite matrix $A_i(\uU_i)$ given by \eqref{defAUi}.
For simplicity of exposition we give the complexity analysis for $d=3$.  To compute $A_1(\uU_1)$ we need first to compute the vectors
$\cT \times (\uu_{j_2,2}\otimes \uu_{j_3,3})$ for $j_2\in[r_2]$ and $j_3\in[r_3]$.  
Each computation of such a vector has complexity $\cO(N)$,
where $N=n_1n_2n_3$.  The number of such vectors is $r_2r_3$.  To form the matrix $A_i(\uU_i)$ we need $\cO(r_2r_3n_1^2)$ flops.
To find the first $r_1$ eigenvectors of $A_1(\uU_1)$ we need $\cO(r_1n_1^2)$ flops.  Assuming that $n_1,n_2,n_3\approx n$ and $r_1,r_2,r_3\approx r$
we deduce that we need $O(r^2 n^3)$  flops to find the first $r_1$ orthonormal eigenvectors of $A_1(\uU_1)$ which span $\U_1(\uU_1)$.
Hence the complexity of finding orthonormal bases of $\U_1(\uU_1)$ is $\cO(r^2n^3)$, which is the complexity of computing $A_1(\uU_1)$.  
Hence the complexity of each step of AMM for best $\br$-approximation,  i.e. computing $\psi^{(l)}$, is $\cO(r^2n^3)$.

It is possible to reduce the complexity of AMM for best $\br$-approximation is to $\cO(r n^3)$ if we compute and store the matrices 
$\cT \times \uu_{j_1,1},\cT \times \uu_{j_2,2},\cT \times \uu_{j_3,3}$.  See \S\ref{sub:compNewtrap}.

We now analyze the complexity of AMM for rank one approximation.  In this case we need only to compute the vector of the form 
$\vv_i:=\cT\times(\otimes_{j\in[d]\setminus\{i\}}\uu_j)$ for each $i\in[d]$, where $\U_j=\span(\uu_j)$ for $j\in[i]$.  The computation of each
$\vv_i$ needs $\cO((d-2)N)$ flops, where $N=\prod_{j\in[d]}$.  Hence each step of AMM for best rank one approximation is $\cO(d(d-2)N)$.
So for $d=3$ and $n_1\approx n_2\approx n_3$ the complexity is $\cO(n^3)$, which is the same complexity as above with $r=1$.
\section{Fixed points of AMM and Newton method}\label{S:fixpt}
Consider the AMM as described in \S\ref{sub:gendp}.  Assume that the sequence $\psi^{(l)}, l\in\N$ converges to a point $\phi\in\Psi$.
Then $\phi$ is a fixed point of the map:
\begin{equation}\label{defFPsi}
\tilde\bF:\Psi\to \Psi, \quad \tilde\bF=(\tilde F_1,\ldots,\tilde F_d),\; \tilde F_i:\Psi\to \Psi_i,\; \tilde F_i(\psi)=\psi_i^\star(\hat\psi_i),\; \psi=(\psi_i,\hat\psi_i),\; i\in[d].
\end{equation}
In general, the map $\tilde \bF$ is a multivalued map, since the maximum given in \eqref{maxprobPsii} may be achieved at a number of points
denoted by $\psi_i^\star(\hat\psi_i)$.  In what follows we assume:
\begin{assum}\label{fixpass}  The AMM converges to a fixed point $\phi$ of $\tilde \bF$ i.e. $\tilde\bF(\phi)=\phi$, such that the following conditions hold:
\begin{enumerate}
\item There is a connected open neighborhood $O\subset \Psi$ such that $\tilde\bF:O\to O$ is one valued map.
\item $\tilde\bF$  is a contraction on $O$ with respect to some norm on $O$.
\item $\tilde\bF\in \rC^2(O)$, i.e. $\tilde\bF$ has two continuous partial derivatives in $O$.
\item $O$ is diffeomorphic to an open subset in $\R^L$.  That is, there exists a smooth one-to-one map $H:O\to \R^L$ such that
the Jacobian $D(H)$ is invertible at each point $\psi\in O$.
\end{enumerate}
\end{assum}

Assume that the conditions of Assumption \ref{fixpass}  hold.   Then the map $\tilde \bF: O\to O$ can be represented
as
\[\bF:O_1\to O_1, \quad \bF=H\circ \tilde \bF \circ H^{-1}, \quad O_1=H(O).\]
Hence to find a fixed point of $\tilde \bF$ in $O$ it is enough to find a fixed point of $\bF$ in $O_1$.
A fixed point of $\bF$ is a zero point of the system
\begin{equation}\label{zeroGx}
\bG(\x)=\0, \quad G(\x):=\x-\bF(\x).
\end{equation}
To find a zero of $\bG$ we use the standard Newton method.
  
In this paper we propose new Newton methods.   We make a few iterations of AMM and switch to a Newton method
assuming that the conditions of Assumption \ref{fixpass} hold as explained above. 
A fixed point of the map $\tilde\bF$ for best rank one approximation induces a fixed point of map $\bF: \R^\n\to \R^\n$ \cite[Lemma 2]{FMPS13}.
Then the corresponding Newton method to find a zero of $\bG$ is straightforward to state and implement, as explained in the next subsection. 
This Newton method was given in \cite[\S5]{FMPS11}.  See also Zhang-Golub \cite{ZhaG01} for a different Newton method for best $(1,1,1)$ approximation.

Let $\tilde\bF:\Gr(\br,\n)\to \Gr(\br,\n)$ be the induced map AMM.  Each $\Gr(r,\R^n)$ can be decomposed as a compact manifold to a finite number 
of charts as explained in \S\ref{S:bestbrapprox}.
These charts induce standard charts of $\Gr(\br,\n)$.
After a few AMM iterations we assume that the neighborhood $O$ of the fixed point of $\tilde\bF$ lies in one the charts of $\Gr(\br,\n)$.
We then construct the corresponding map $\bF$ in this chart.  Next we apply the standard Newton method to $\bG$.  
The papers by Eld\'en-Savas \cite{ES09} and Savas-Lim \cite{SL10} discuss Newton and quasi-Newton methods for $(r_1,r_2,r_3)$ approximation
of $3$-tensors using the concepts of differential geometry.

\subsection{Newton method for best rank one approximation}\label{sub:newtboneappr}
Let $\cT\in\R^{\n}\setminus \{0\}$.  Define:
\begin{eqnarray}\notag
&&\Psi_i=\R^{n_i},\; i\in [d],\quad \Psi=\R^{n_1}\times \cdots \times \R^{n_d}, \quad \psi=(\x_1,\ldots,\x_d)\in \Psi,\\
&&f_\cT:\Psi\to \R, \quad f_\cT(\psi)=\cT\times (\otimes_{j\in[d]}\x_j),\\
&&\bF=(F_1,\ldots,F_d):\Psi\to \Psi, \quad F_i(\psi)=\cT\times (\otimes_{j\in[d]\setminus\{i\}}\x_j), \quad i\in [d].\label{defbFfcT}
\end{eqnarray}
Recall the results of \S\ref{sub:singvten}: Any critical point of  $f_\cT|\rS(\n)$ satisfies \eqref{defsingvalstup}.  Suppose we start the AMM with $\psi^{(0)}=(\x_1^{(0)},\ldots,\x^{(0)}_d)\in\rS(\n)$ such that
$f_\cT(\psi^{(0)})\ne 0$.  Then it is straightforward to see that $f_\cT(\psi^{(l)})>0$ for $l\in\N$.  Assume that $\lim_{l\to\infty}\psi^{(l)}=\omega=(\uu_1,\ldots,\uu_d)
\in\rS(\n)$.   Then $\omega$ is the singular tuple of $\cT$ satisfying \eqref{defsingvalstup}.  Clearly, $\lambda=f_\cT(\omega)>0$.
Let  
\begin{equation}\label{defnorphi}
\phi=(\y_1,\ldots,\y_d):=\lambda^{-\frac{1}{d-2}}\omega=\lambda^{-\frac{1}{d-2}}(\uu_1,\ldots,\uu_d).
\end{equation}
Then $\phi$ is a fixed point of $\bF$.  

Our Newton algorithm for finding the fixed point $\phi$ of $\bF$ corresponding to a fixed point $\omega$ of AMM is as follows.  We do a number of iterations of
AMM to obtain $\psi^{(m)}$.  Then we renormalize $\psi^{(m)}$ according to \eqref{defnorphi}:
\begin{equation}\label{phi0for}
\phi_0:=(f_\cT(\psi^{(m)})^{-\frac{1}{d-2}}\psi^{(m)}.
\end{equation}
Let $D\bF(\psi)$ denote the Jacobian of $\bF$ at $\psi$, i.e. the matrix of partial derivatives of $\bF$ at $\psi$.
Then we perform Newton iterations of the form:
\begin{equation}\label{Newtitomeg}
\phi^{(l)}=\phi^{(l-1)} -(I-D\bF(\phi^{(l-1)}))^{-1}(\phi^{(l-1)}-\bF(\phi^{(l-1)})), \quad l\in\N.
\end{equation}

After performing a number of Newton iterations we obtain $\phi^{(m')}=(\z_1,\ldots,\z_d)$ which is an approximation of $\phi$.  We then renormalize 
each $\z_i$ to obtain $\omega^{(m')}:=(\frac{1}{\|\z_1\|}\z_1,\ldots,\frac{1}{\|\z_d\|}\z_d)$ which is an approximation to the fixed point $\omega$.
We call this Newton method \emph{Newton-1}.

We now give the explicit formulas for $3$-tensors, where $n_1=m,n_2=n,n_3=l$.  First
\begin{equation}\label{defmapF}
\bF(\uu,\vv,\w):=(\cT\times(\vv\otimes\w),\cT\times(\uu\otimes\w),\cT\times(\uu\otimes\vv)), \quad \bG:=(\uu,\vv,\w)-\bF(\uu,\vv,\w).
\end{equation}
Then
\begin{equation}\label{JacDF}
 D\bG(\uu,\vv,\w)=\left[\begin{array}{ccc}I_m&-\cT\times \w&-\cT\times\vv\\-(\cT\times\w)\trans&I_n&-\cT\times\uu\\
 -(\cT\times \vv)\trans&-(\cT\times \uu)\trans&I_l\end{array}\right].
\end{equation}
 Hence Newton-1 iteration is given by the formula
\[(\uu_{i+1},\vv_{i+1},\w_{i+1})= (\uu_{i},\vv_{i},\w_i) -(D\bG(\uu_i,\vv_i,\w_i))^{-1}\bG(\uu_i,\vv_i,\w_i),\]
 for $i=0,1,\ldots,$.
Here we abuse notation by viewing $(\uu,\vv,\w)$ as a column vector $(\uu\trans,\vv\trans,\w\trans)\trans\in \C^{m+n+l}$.

Numerically, to find $(D\bG(\uu_i,\vv_i,\w_i))^{-1}\bG(\uu_i,\vv_i,\w_i)$ one solves the linear system
\[
 (D\bG(\uu_i,\vv_i,\w_i))(\x,\y,\z)=\bG(\uu_i,\vv_i,\w_i).
\]
The final vector $(\uu_j,\vv_j,\w_j)$ of Newton-1 iterations is
followed by a scaling to vectors of unit length $\x_j=\frac{1}{\|\uu_j\|}\uu_j, \y_j=\frac{1}{\|\vv_j\|}\vv_j, \z_j=\frac{1}{\|\w_j\|}\w_j$.

We now discuss the complexity of Newton-1 method for $d=3$.  Assuming that $m\approx n\approx l$ we deduce that the computation of
the matrix $D\bG$ is $\cO(n^3)$. As the dimension of $D\bG$ is $m+n+l$ it follows that the complexity of each iteration of Newton-1 method is
$\cO(n^3)$.
\section{Newton method for best $\br$-approximation}\label{S:bestbrapprox}
Recall that an $r$-dimensional subspace $\U\in \Gr(r,\R^n)$ is given by a matrix $U=[u_{ij}]_{i,j=1}^{n,r}\in\R^{n\times r}$
of rank $r$.  In particular there is a subset $\alpha\subset [n]$ of cardinality $r$ so that $\det U[\alpha,[r]]\ne 0$.  Here $\alpha=(\alpha_1,\ldots,\alpha_r),
1\le \alpha_1<\ldots <\alpha_r\le n$.  So $U[\alpha,[r]]:=[u_{\alpha_i j}]_{i,j=1}^r\in \gl(r,\R)$, (the group of invertible matrices).
Clearly, $V:=U U[\alpha,[r]]^{-1}$ represents another basis in $\U$.  Note that $V[\alpha,[r]]=I_r$.   
Hence the set of all $V\in \R^{n\times r}$ with the condition: $V[\alpha,[r]]=I_r$ represent an open cell in $\Gr(r,n)$ of dimension $r(n-r)$ 
denoted by $\Gr(r,\R^n)(\alpha)$.  
(The number of free parameters in all such $V$'s is $(n-r)r$.)  Assume for simplicity of exposition that $\alpha=[r]$.  Note that 
$V_0=\left[\begin{array}{c}I_r\\0\end{array}\right]
\in \Gr(r,\R^n)([r])$.  Let $\e_i=(\delta_{1i},\ldots,\delta_{ni})\trans\in \R^n, i=1,\ldots,n$ be the standard basis in $\R^n$.  So 
$\U_0=\span(\e_1,\ldots,\e_r)\in \Gr(r,\R^n)([r])$,
and $V_0$ is the unique representative of $\U_0$.  Note that  $\U_0^\perp$, the orthogonal complement of $\U_0$, is $\span(\e_{r+1},\ldots,\e_n)$. 
It is straightforward to see that $\V\in \Gr(r,\R^n)([r])$ if and only if $\V\cap \U_0^\perp=\{\0\}$.

 The following definition is a geometric generalization of $\Gr(r,\R^n)(\alpha)$:
\begin{equation}\label{openUcell}
\Gr(r,\R^n)(\U):=\{\V\in \Gr(r,\R^n),\; \V\cap \U^\perp=\{\0\}\} \textrm{ for } \U\in\Gr(r,\R^n).
\end{equation}
A basis for $\Gr(r,\R^n)(\U)$, which can be identified the tangent hyperplane $T_\U \Gr(r,\R^n)$, can be represented as $\oplus^r \U^\perp$:
Let $\uu_1,\ldots,\uu_r$ and $\uu_{r+1},\ldots,\uu_n$ be orthonormal bases  of $\U$ and $\U^\perp$ respectively
Then each subspace $\V\in \Gr(r,\R^n)(\U)$ has a unique basis of the form $\uu_1+\x_1,\ldots,\uu_r+\x_r$ for unique $\x_1,\ldots,\x_r\in \U^\perp$.
Equivalently, every matrix $X\in \R^{(n-r)\times r}$ induces a unique subspace $\V$ using the equality
\begin{equation}\label{defx1xr}
[\x_1\;\ldots\;\x_r]=[\uu_1\;\ldots\;\uu_{n-r}] X \textrm{ for each } X\in\R^{(n-r)\times r}.
\end{equation}

 Recall the results of \S\ref{sub:apprprb}.
 Let $\underline{U}=(U_1,\ldots,U_d)\in\Gr(\br,\n)$.  Then 
\begin{equation}\label{deftilFGrn}
\tilde\bF=(\tilde F_1,\ldots, \tilde F_d):\Gr(\br,\n)\to\Gr(\br,\n), \quad \tilde F_i(\uU)=\U_i(\uU_i), \; i\in[d],
\end{equation}
where $\U_i(\uU_i)$ a subspace spanned by the first $r_i$ eigenvectors of $A_i(\uU_i)$.
Assume that $\tilde\bF$ is one valued at  $\uU$, i.e. 
\eqref{uniqcondUiAi} holds.  Then it is straightforward to show that $\tilde\bF$ is smooth (real analytic) in neighborhood of $\underline{U}$.
Assume next that there exists a neighborhood $O$ of $\underline{U}$ such that 
\begin{equation}\label{regFcond}
 O\subset  \Gr(\br,\n)(\uU):=\Gr(r_1,\R^{n_1})(\U_1)\times \cdots \times\Gr(r_d,\R^{n_d})(\U_d), \quad \uU=(\U_1,\ldots,\U_d),
\end{equation}
such that the conditions \emph{1-3} of Assumption \ref{fixpass} hold.
Observe next that $\Gr(\br,\n)(\uU)$ is diffeomorphic to 
\[\R^L:= \R^{(n_1-r_1)\times r_1} \ldots \times \R^{(n_d-r_d)\times r_d}, \quad L=\sum_{i\in[d]} (n_i-r_i)r_i\]
We say that $\tilde\bF$ is \emph{regular} at $\underline{U}$ if in addition to the above condition the matrix $I-D\tilde \bF(\uU)$ is invertible.
We can view $X=[X_1\;\ldots \;X_d]\in  \R^{(n_1-r_1)\times r_1} \ldots \times \R^{(n_d-r_d)\times r_d}$.
Then $\tilde \bF$ on $O$ can be viewed as
\begin{equation}\label{tildebFre}
\bF:O_1\to O_1, \quad O_1\subset \R^L,\quad \bF(X)=[F_1(X),\ldots, F_d(X)],\; X=[X_1\;\ldots \;X_d]\in\R^L.
\end{equation}
Note that $F_i(X)$ does not depend on $X_i$ for each $i\in[d]$.  
In our numerical simulations we first do a small number of AMM and then switch to Newton method given by \eqref{Newtitomeg}.
Observe that  $\uU$ corresponds to $X(\uU)=[X_1(\uU),\ldots,X_d(\uU)]$.  
When referring to \eqref{Newtitomeg} we identify $X=[X_1,\ldots,X_d]$ with $\phi=(\phi_1,\ldots,\phi_d)$ and no ambiguity will arise.

Note that the case $\br=\1_d$ corresponds to best rank one approximation.  The above Newton method in this case is different from Newton 
method given in \S\ref{sub:newtboneappr}.
\section{A closed formula for $D\bF(X(\uU))$}\label{S:forDF}
Recall the definitions and results of \S\ref{sub:apprprb}.  Given $\uU$ we compute $\tilde F_i(\uU)=\U_i(\uU_i)$, which is the subspace
spanned by the first $r_i$ eigenvectors of $A_i(\uU_i)$, which is given by \eqref{defAUi}, for $i\in[d]$.  Assume that \eqref{basisUi} holds.
Let
\begin{eqnarray}\notag
&&\U_{i}(\uU_i)=\span(\vv_{1,i},\ldots,\vv_{r_i,n_i}), \quad \U_i(\uU_i)^\perp=\span(\vv_{r_i+1,i},\ldots, \vv_{n_i,i}),\\
&&\vv_{j,i}\trans \vv_{k,i}=\delta_{jk},\;j,k\in [n_i], \quad i\in [d].\label{basistFUi}
\end{eqnarray}

With each $X=[X_1,\ldots,X_d]\in\R^L$ we associate the following point $(\W_1,\ldots,\W_d)\in\Gr(\br,\n)(\uU)$.
Suppose that $X_i=[x_{pq,i}]\in\R^{(n_i-r_i)\times r_i}$.  Then $\W_i$ has a basis of the form
\[\uu_{j_i,i}+\sum_{k_i\in [n_i-r_i]}x_{k_ij_i,i}\uu_{r_i+k_i,i}, \quad j_i\in [r_i].\]
One can use the following notation for a basis $\w_{1,i},\ldots,\w_{r_i,i}$, written as a vector with vector coordinates $[\w_{1,i}\cdots\w_{r_i,i}]$:
\begin{equation}\label{basisWi}
[\w_{1,i}\cdots\w_{r_i,i}]=[\uu_{1,i}\cdots\uu_{r_i,i}] + [\uu_{r_i+1,i}\cdots\uu_{n_i,i}]X_i, \quad i\in[d].
\end{equation}

Note that to the point $\uU\in \Gr(\br,\n)(\uU)$ corresponds the point $X=0$.  Since $\uu_{1,i},\ldots,\uu_{n_i,i}$ is a basis in $\R^{n_i}$
it follows that
\begin{eqnarray}\notag
&&[\vv_{1,i}\cdots \vv_{r_i,i}]=[\uu_{1,i}\cdots \uu_{r_i,i}]Y_{i,0}+[\uu_{r_i+1,i},\ldots,\uu_{n_i,i}]X_{i,0}=[\uu_{1,i}\cdots \uu_{n_i,i}]Z_{i,0}, \\
&& Y_{i,0}\in \R^{r_i\times r_i}, \;X_{i,0}\in \R^{(n_i-r_i)\times r_i},\;Z_{i,0}=\left[\begin{array}{c}Y_{i,0}\\X_{i,0}\end{array}\right]\in \R^{n_i\times r_i},
\textrm{ for }i\in [d].\label{vforuXY}
\end{eqnarray}
View $[\vv_{1,i}\cdots \vv_{r_i,i}],[\uu_{1,i}\cdots \uu_{n_i,i}]$ as $n_i\times r_i$ and $n_i\times n_i$ matrices with orthonormal columns.  Then
\begin{equation}\label{forZi0}
Z_{i,0}=[\uu_{1,i}\cdots \uu_{n_i,i}]\trans[\vv_{1,i}\cdots \vv_{r_i,i}], \quad i\in [d].
\end{equation}

The assumption that $\tilde \bF:O\to O$ implies that $Y_{i,0}$ is an invertible matrix.  Hence $[\vv_{1,i}\cdots \vv_{r_i,i}]Y_{i,0}^{-1}$ is also a basis 
in $\U_i(\uU_i)$.  Clearly,
\[[\vv_{1,i}\cdots \vv_{r_i,i}]Y_{i,0}^{-1}=[\uu_{1,i}\cdots \uu_{r_i,n_i}]+[\uu_{r_i+1,i},\ldots,\uu_{n_i,i}]X_{i,0}Y_{i,0}^{-1}, \quad i\in [d].\]

Hence $\tilde\bF(\uU)$ corresponds to $\bF(0)$ where
\begin{equation}\label{exprebF0}
 F_i(0)=X_{i,0}Y_{i,0}^{-1}, \;i\in [d], \quad \bF(0)=(F_1(0),\ldots,F_d(0)).
\end{equation}

We now find the matrix of derivatives.  So $D_iF_j\in\R^{((n_i-r_i)r_i\times ((n_j-r_j)r_j}$ is the partial derivative matrix of $(n_j-r_j)r_j$ coordinates
of $F_j$ with respect to $(n_i-r_i)r_i$ the coordinates of $\U_i$ viewed as the matrix $\left[\begin{array}{c}I_{r_i}\\G_i\end{array}\right]$.
So $G_i\in \R^{(n_i-r_i)\times r_i}$ are the variables representing the subspace $\U_i$.  Observe first that
$D_iF_i=0$ just as in Newton method for best rank one approximation in \S\ref{sub:newtboneappr}.  

Let us now find $D_i F_j(0)$.  Recall that $D_i F_j(0)$ is a matrix of size $(n_i-r_i)r_i\times (n_j-r_j)r_j$.
The entries of  $D_i F_j(0)$ are indexed by $((p,q),(s,t))$ as follows: The entries of $G_i=[g_{pq,i}]\in \R^{(n_i-r_i)\times r_i}$
are viewed as $(n_i-r_i)r_i$ variables, and are indexed by $(p,q)$, where $p\in[n_i-r_i],q\in[r_i]$.
$F_j$ is viewed as a matrix $G_j\in \R^{(n_j-r_j)\times r_j}$.The entries of $F_j$ are indexed by $(s,t)$, where $s\in[n_j-r_j]$ and $t\in[r_j]$. 
Since $\uU\in\Gr(\br,\n)(\uU)$ corresponds to $0\in\R^L$ we denote by $A_j(0)$ the matrix $A_j(\uU_j)$ for $j\in[d]$.
We now give the formula for $\frac{\partial  A_j(0)}{\partial g_{pq,i}}$.  This is done by noting that we vary $\U_{i}$ by changing the orthonormal
basis $\uu_{1,i},\ldots,\uu_{r_i,i}$ up to the first perturbation with respect to the real variable $\varepsilon$ to
\[\hat\uu_{1,i}=\uu_{1,i},\ldots,\hat \uu_{q-1,i}= \uu_{q-1,i},\hat\uu_{q,i}=\uu_{q,i}+\varepsilon\uu_{r_i+p,i},
\hat\uu_{q+1,i}=\uu_{q+1,i},\ldots,\hat\uu_{r_i,i}=\uu_{r_i,i}\]
We denote the subspace spanned by these vectors as $\U_{i}(\varepsilon,p,q)$.
That is, we change only the $q$ orthonormal vector of the standard basis in $\U_{i}$, for $q=1,\ldots,r_i$.  The new basis is an orthogonal basis, 
and up order $\varepsilon$, 
the vector $ \uu_{q,i}+\varepsilon\uu_{r_i+p,i}$ is also of length $1$.  
Let $\uU(\varepsilon,i,p,q)=(\U_1,\ldots,\U_{i-1},\U_i(\varepsilon,p,q),\U_{i+1},\ldots,\U_d)$. 
Then $\uU(\varepsilon,i,p,q)_j$ is obtained by dropping the subspace $\U_j$ from  $\uU(\varepsilon,i,p,q)$.
We will show that
\begin{equation}\label{defAjvareps}
A_j(\uU(\varepsilon,i,p,q)_j)=A_j(\uU_j)+\varepsilon B_{j,i,p,q}+O(\varepsilon^2).
\end{equation}

We now give a formula to compute $B_{j,i,p,q}$.
Assume that  $i,j\in [d]$ is a pair of different integers.
Let $J$ be a set of $d-2$ pairs $\cup_{l\in [d]\setminus\{i,j\}}\{(k_l,l)\}$, where $k_l\in[r_l]$.
Denote by $\cJ_{ij}$ the set of all such $J$'s.  Note that $\cJ_{ij}=\cJ_{ji}$.
Furthermore, the number of elements in $\cJ_{ij}$ is $R_{ij}=\prod_{l\in[d]\setminus\{i,j\}} r_l$.
We now introduce the following matrices
\begin{equation}\label{defCijk}
C_{ij}(J):=\cT\times (\otimes_{(k,l)\in J}\uu_{k,l})\in \R^{n_i\times n_j}, \quad J\in\cJ_{ij}.
\end{equation}
Note that $C_{ij}(J)=C_{ji}(J)\trans$.  
\begin{lemma}\label{AjBjpqfor}  Let $i,j\in [d],i\ne j$.  Assume that $p\in[n_i-r_i],q\in[r_i]$.  Then \eqref{defAjvareps} holds.
Furthermore
\begin{eqnarray}\label{newforAj}
&&A_j(\uU_j)=\sum_{k\in[r_i], J\in\cJ_{ji}} (C_{ji}(J)\uu_{k,i})(C_{ji}(J)\uu_{k,i})\trans,\\
&&B_{j,i,p,q}=\sum_{J\in\cJ_{ji}} (C_{ji}(J)\uu_{k_i+p,i})(C_{ji}(J)\uu_{q,i})\trans+(C_{ji}(J)\uu_{q,i})(C_{ji}(J)\uu_{k_i+p,i})\trans
\label{forBjipq}
\end{eqnarray}
\end{lemma}
\proof  The identity of \eqref{newforAj} is just a restatement of \eqref{defAUi}.  To compute $A_j(\uU(\varepsilon,i,p,q)_j)$ use \eqref{defCijk} by replacing 
$u_{k,i}$ with $\hat u_{k,i}$ for $k\in[n_i]$.  Deduce first \eqref{defAjvareps} and then \eqref{forBjipq}.\qed
 
Recall that $\vv_{1,j},\ldots,\vv_{r_j,j}$ is an orthonormal basis of $\U_j(\uU_j)$, and these vectors are the eigenvectors $A_j(\uU_j)$ corresponding 
its first $r_j$ eigenvalues.   Let $\vv_{r_j+1,j},\ldots,\vv_{n_j,i}$ be the last $n_j-r_j$ orthonormal eigenvectors of $A_j(\uU_j)$.
We now find the first perturbation of the first $r_i$ eigenvectors for the matrix $A_j(\uU_j)+\varepsilon B_{j,i,p,q}$.
Assume first, for simplicity of exposition, that each $\lambda_k(A_j(\uU_j))$ is simple for $k\in [r_j]$:
Then it is known, e.g. \cite[Chapter 4, \S19, (4.19.2)]{Fri15}:
\begin{equation}\label{pertforxj2}
\vv_{k,j}(\varepsilon,i,p,q)=\vv_{k,j}+\varepsilon (\lambda_{k}(A_j(\uU_j)) I_{n_j}-A_j(\uU_j))^\dagger B_{j,i,p,q}\vv_{k,j} +O(\varepsilon^2),  k\in[r_j].
\end{equation}

The assumption that  $\lambda_{k}(A_j(\uU_j))$ is a simple eigenvalue for $k\in[r_j]$ yields
\[(\lambda_{k}(A_j(\uU_i))I_{n_j}-A_j(\uU_j))^\dagger\y=\sum_{l\in[n_j] 
\setminus\{k\}}\frac{1}{\lambda_{k}(A_j(\uU_j))-\lambda_{l}(A_j(\uU_j))} (\vv_{l,j}\trans\y)\vv_{l,j},\]
for $\y\in\R^{n_j}$.

Since we are interested in a basis of $\U_j(\uU(\varepsilon,i,p,q)_j)$ up to the order of $\varepsilon$ we can assume that this basis is of the form
\[\tilde \vv_{k,j}(\varepsilon,i,p,q)=\vv_{k,j}+\varepsilon\w_{k,j}(i,p,q), \quad \w_{k,j}(i,p,q)\in
\span(\vv_{r_j+1,j}\ldots,\vv_{n_j,j}).\]
Hence
\begin{eqnarray}\notag
 &&\w_{k,j}(i,p,q)=\sum_{l\in[n_j]\setminus [r_j]}\frac{1}{\lambda_{k}(A_j(\uU_j))-\lambda_{l}(A_j(\uU_j))} (\vv_{l,j}\trans \c_{k,j,i,p,q})\vv_{l,j},\\
&&\c_{k,j,i,p,q}:= B_{j,i,p,q}\vv_{k,j}.\label{pertforbasU2}
\end{eqnarray}
Note that the assumption \eqref{uniqcondUiAi} yields that $\w_{k,j}$ is well defined for $k\in[r_j]$.  Let
\begin{eqnarray*}
&&W_j(i,p,q)=[\w_{1,j}(i,p,q)\cdots\w_{r_j,j}(i,p,q)]=\left[\begin{array}{c}V_j(i,p,q)\\U_j(i,p,q)\end{array}\right],\\ 
&&V_j(i,p,q)\in\R^{r_j\times r_j},\quad U_j(i,p,q)\in\R^{(n_j-r_j)\times r_j}.
\end{eqnarray*}
Up to the order of $\varepsilon$ we have that a basis of $\U_{j}(\uU(\varepsilon,i,p,q)_j)$ is given by columns of matrix 
$Z_{j,0}+\varepsilon W_j(i,p,q)=\left[\begin{array}{c}Y_{j,0}+\varepsilon V_j(i,p,q)\\X_{j,0} +\varepsilon U_j(i,p,q)\end{array}\right]$.  Note
\[(Z_{j,0}+\varepsilon W_j(i,p,q))(Y_{j,0}+\varepsilon V_j(i,p,q))^{-1}=\left[\begin{array}{c} I_{r_j}\\
(X_{j,0}+\varepsilon U_j(i,p,q))(Y_{j,0}+\varepsilon V_j(i,p,q))^{-1}\end{array}\right].\]
Observe next
\begin{eqnarray*}
&&Y_{j,0}+\varepsilon V_j(i,p,q)=Y_{j,0}(I_{r_j}+\varepsilon Y_{j,0}^{-1} V_j(i,p,q)), \\
&&(Y_{j,0}+\varepsilon V_j(i,p,q))^{-1}=(I_{r_j}+\varepsilon Y_{j,0}^{-1} V_j(i,p,q))^{-1}Y_{j,0}^{-1}=\\
&&Y_{j,0}^{-1}-\varepsilon Y_{j,0}^{-1} V_j(i,p,q)Y_{j,0}^{-1}+O(\varepsilon^2),\\
&&(X_{j,0}+\varepsilon U_j(i,p,q))(Y_{j,0}+\varepsilon V_j(i,p,q))^{-1}=\\
&&X_{j,0}Y_{j,0}^{-1}+\varepsilon(U_j(i,p,q) Y_{j,0}^{-1}-X_{j,0}Y_{j,0}^{-1}V_{j}(i,p,q) Y_{j,0}^{-1})+O(\varepsilon^2).
\end{eqnarray*} 
Hence
\begin{equation}\label{D1F2pqfor}
\frac{\partial F_j}{\partial g_{pq,i}}(0)=U_j(i,p,q)Y_{j,0}^{-1}-X_{j,0}Y_{j,0}^{-1}V_{j}(i,p,q) Y_{j,0}^{-1}.
\end{equation}

Thus $D\bF(0)=[D_iF_J]_{i,j\in[d]}\in\R^{L\times L}$.  We now make one iteration of Newton method given by \eqref{Newtitomeg}
for $l=1$, where $\phi^{(0)}=0$:
\begin{equation}\label{Newtmetstrt0}
\phi^{(1)}=-(I-D\bF(0))^{-1}F(0), \quad \phi^{(1)}=[X_{1,1},\ldots,X_{d,1}]\in\R^L.
\end{equation}
Let $\U_{i,1}\in\Gr(r_i,\R^{n_i})$ be the subspace represented by the matrix $X_{i,1}$: 
\begin{equation}\label{defUi1}
\U_{i,1}=\span(\tilde\uu_{1,i,1},\ldots,\tilde\uu_{r_i,i,1}),\; [\tilde\uu_{1,i,1},\ldots,\tilde\uu_{n_i,i,1}]=
[\uu_{1,i},\ldots,\uu_{n_i,i}]\left[\begin{array}{c} I_{r_i}\\X_{j,1}\end{array}\right]
\end{equation}
for $i\in[d]$.  Perform the Gram-Schmidt process on $\tilde\uu_{1,i,1},\ldots,\tilde\uu_{r_i,i,1}$ to obtain an orthonormal basis
$\uu_{1,i,1},\ldots,\uu_{r_i,i,1}$ of $\U_{i,1}$.  Let $\uU:=(\U_{1,1},\ldots,\U_{d,1})$ and repeat the algorithm which is described above.
We call this Newton method \emph{Newton-2}.
\section{Complexity of Newton-2}\label{sub:compNewtrap}
In this section we assume for simplicity that $d=3$, $r_1=r_2=r_3=r$, $n_i\approx n$ for $i\in [3]$.
We assume that executed a number of times the AMM for a given $\cT\in \R^\n$.  So we are given $\uU=(\U_1,\ldots,\U_d)$, and an orthonormal basis
$\uu_{1,i},\ldots,\uu_{r,i}$ of $\U_i$ for $i\in[d]$.  
First we complete each $\uu_{1,i},\ldots,\uu_{r,i}$ to an orthonormal basis $\uu_{1,i},\ldots,\uu_{n_i,i}$of $\R^{n_i}$, which needs $\cO(n^3)$ flops.
Since $d=3$ we still need only $\cO(n^3)$ to carry out this completion for each $i\in [3]$.

Next we compute the matrices $C_{ij}(J)$.  Since $d=3$, we need $n$ flops to compute each entry
of $C_{ij}(J)$.  Since we have roughly $n^2$ entries, the complexity of computing  $C_{ij}(J)$ is $\cO(n^3)$.  As the cardinality of $\cJ_{ij}$ is $r$ we need 
$\cO(rn^3)$ flops to compute all $C_{ij}(J)$ for $J\in\cJ_{ij}$.  As the number of pairs in $[3]$ is $3$ it follows that the complexity of computing all $C_{ij}(J)$
is $\cO(rn^3)$.  

The identity  \eqref{newforAj} yields that the complexity of computing $A_j(\uU_j)$ is $\cO(r^2n^2)$.  Recall next that $A_j(\uU_j)$ is $n_j\times n_j$
symmetric positive semi-definite matrix.  The complexity of computations of  the eigenvalues and the orthonormal eigenvectors of $A_j(\uU_j)$ is $\cO(n^3)$.
Hence the complexity of computing $\uU$ is $\cO(r n^3)$, as we pointed out at the end of \S\ref{sub:companbrap}.

The complexity of computing $B_{j,i,p.q}$ using \eqref{forBjipq} is $\cO(r n^2)$.   The complexity of computing  $\w_{k,j}(i,p,q)$, given by \eqref{pertforbasU2}
is $\cO(n^2)$.  Hence the complexity of computing $W_j(i,p,q)$ is $\cO(rn^2)$.  Therefore the complexity of computing $D_i F_j$ is $\cO(r^2n^3)$.
Since $d=3$, the complexity of computing the matrix $D\bF(0)$ is also  $\cO(r^2n^3)$.

As $D\bF(0)\in \R^{L\times L}$, where $L\approx 3r n$, the complexity of computing $(I-D\bF(0))^{-1}$ is $\cO(r^3n^3)$.
In summary, the complexity of one step in Newton-2 is $\cO(r^3 n^3)$.

\section{Numerical Results}\label{S:NumRes}

We have implemented a Matlab library tensor decomposition using Tensor Toolbox given by \cite{KB09}. The performance was 
measured via the actual CPU-time (seconds) needed to compute.
   All performance tests have been carried out on a 2.8 GHz Quad-Core Intel Xeon Macintosh computer with 16GB RAM. The performance results are discussed for 
real data sets of third-order tensors. We worked with a real computer tomography (CT) data set (the so-called MELANIX data set of OsiriX) \cite{FMPS13}.

Our simulation results are averaged over 10 different runs of the each algorithm.  In each run, we changed the initial guess, that is, we generated new random 
start vectors.  We always initialized the algorithms by random start vectors, because this is cheaper 
than the initialization via HOSVD. 
We note here that for Newton methods our initial guess is the subspaces returned by one iteration of AMM method.

All the alternating algorithms have the same stopping criterion where convergence is achieved if one of the two following conditions are met: 
$iterations > 10;fitchange <0.0001$ is met. 
All the Newton algorithms have the same stopping criterion  where convergence is achieved if one of the two following conditions are met: $iterations > 10;change <\exp(-10)$.

Our numerical simulations demonstrate the well known fact that for large size tensors Newton methods are not efficient. Though the Newton methods converge in fewer iterations than alternating methods, the computation associated with the matrix of derivatives (Jacobian) in each iteration is too expensive making alternating maximization methods much more cost effective. 
Our simulations also demonstrate that our Newton-1 for best rank one approximation is as fast as AMM methods.  However our Newton-2  is much 
slower than alternating methods. We also give a comparison between our Newton-2 
and the Newton method based on Grassmannian manifold  by \cite{ES09}, abbreviated as Newton-ES. 

We also observe that  for large tensors and large rank approximation two alternating maximization methods, namely MAMM and 2AMMV, seem to outperform the other alternating maximization methods.
We would recommend Newton-1 for rank one approximation in case of rank one approximation both for large and small sized tensors.  For higher rank approximation we recommend 2AMMV  in case of large size tensors and AMM or MAMM in case of small size tensors.

Our Newton-2 performs a bit slower than Newton-ES,  however we would like to point couple of advantages. Our method can be easily extendable to higher 
dimensions ( for $d>3$ case) both analytically and numerically compared to Newton-ES.
Our method is also highly parallelizable which can  bring down the computation time drastically. Computation of $D_{i}F_{j}$ matrices  in each iteration contributes 
to about $50\%$ of the total time, which however can be parallelizable.   Finally the number of iterations in Newton-2 is at least $30\%$ less than in Newton-ES.

It is not only important to check how fast the different algorithms perform but also what quality
they achieve. This was measured by checking the Hilbert-Schmidt norm, abbreviated as HS norm, of the resulting decompositions,
which serves as a measure for the quality of the approximation. In general, we can say that the higher
the HS norm, the more likely it is that we find a global maximum. Accordingly, we compared
the HS norms to say whether the different algorithms converged to the same stationary point.
In Figure 4, we show the average HS norms achieved by different algorithms and compared them with the input norm. 
We observe all the algorithms seem to attain the same local maximum.

 \includegraphics[scale = .9] {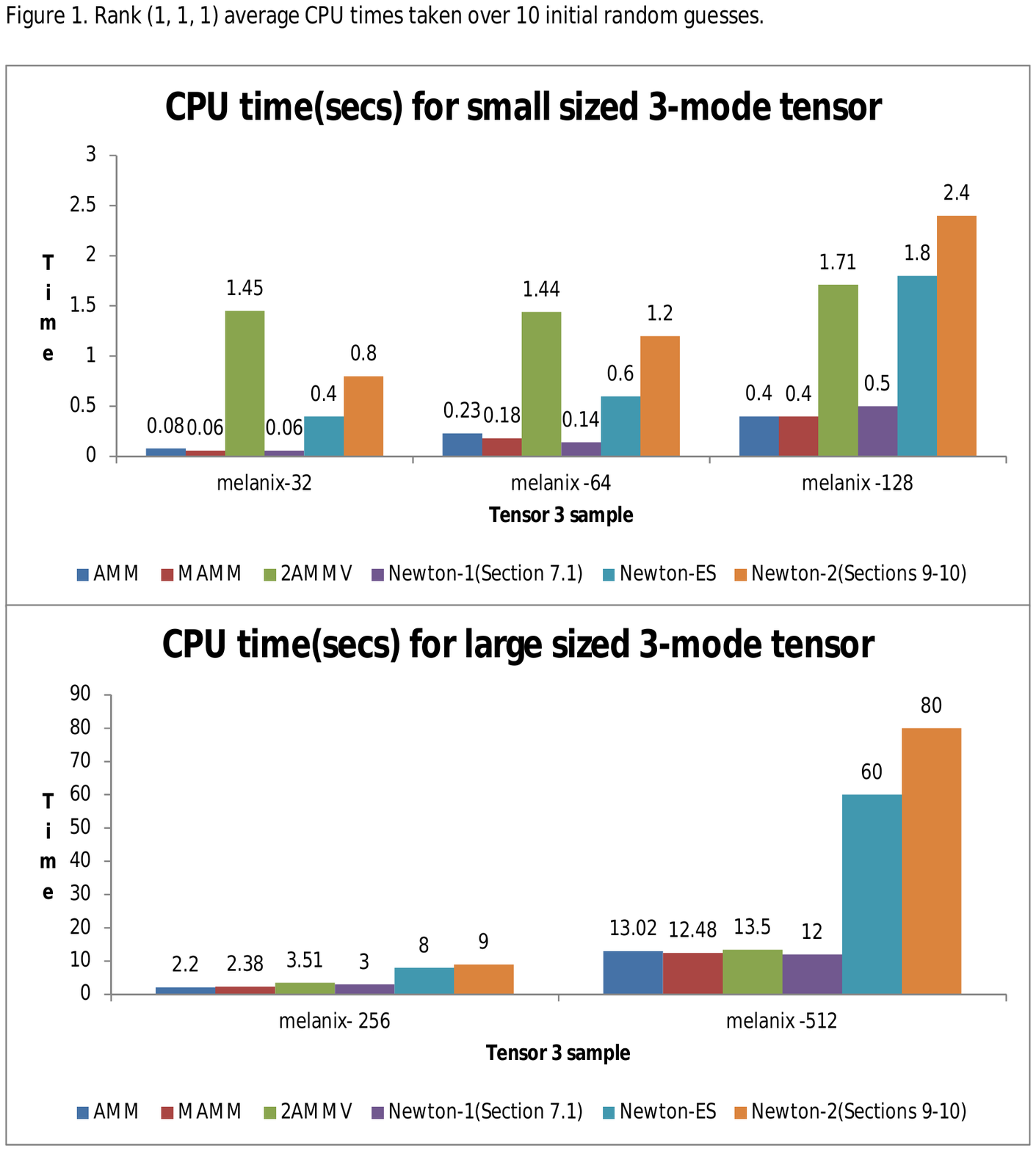}
 \includegraphics[scale = .9] {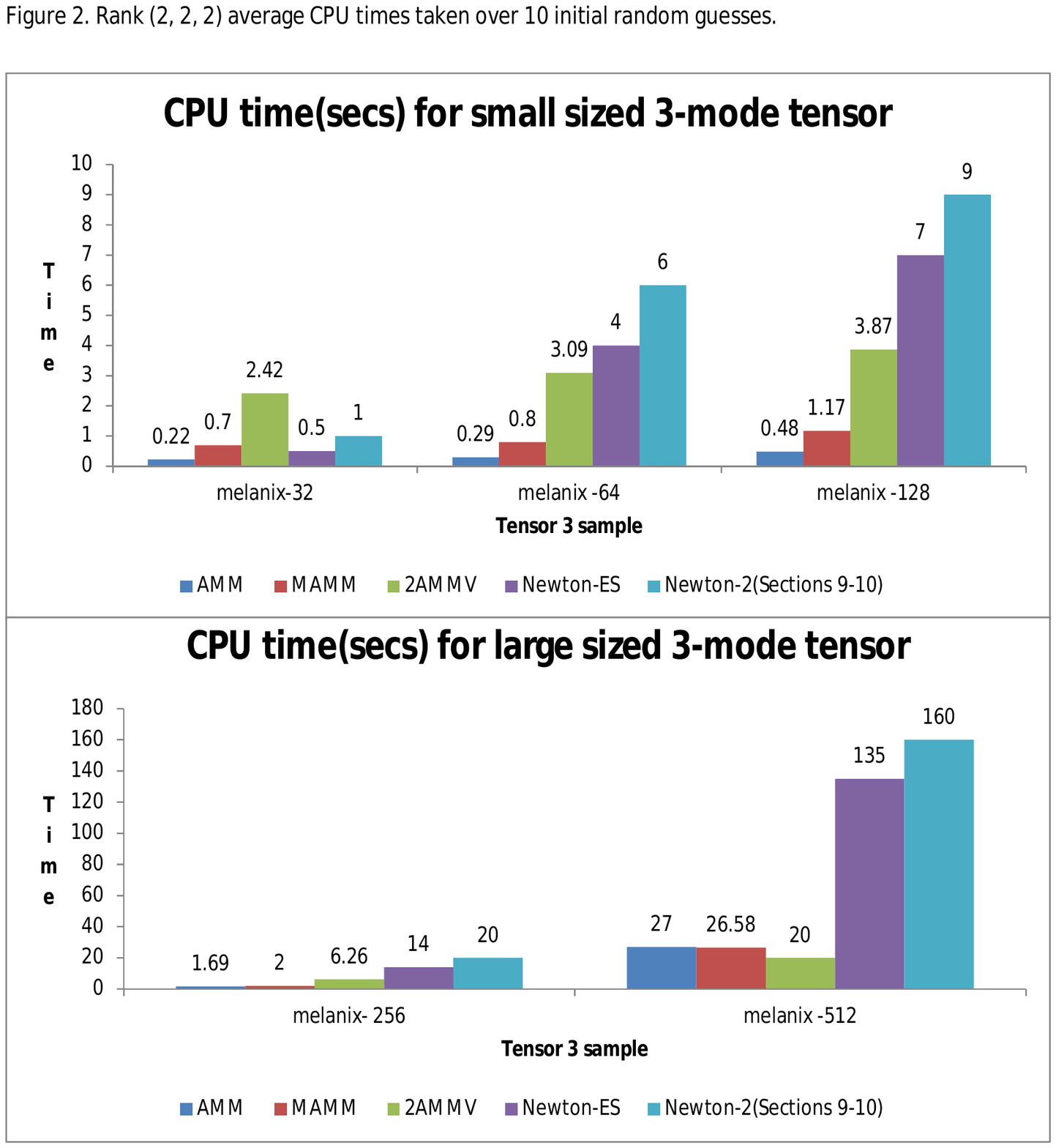}
 \includegraphics[scale = .9] {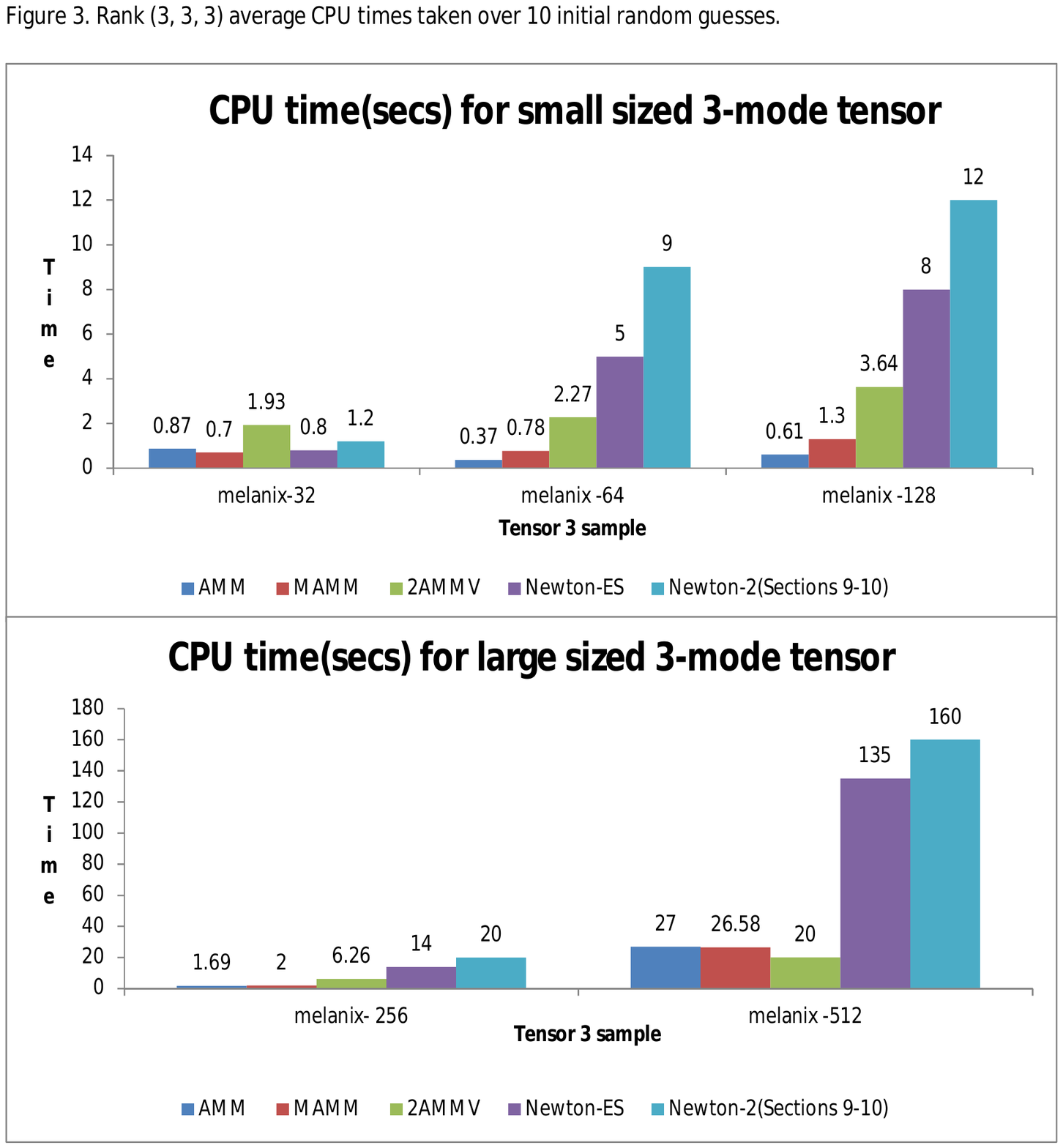}
 \includegraphics[scale = .9] {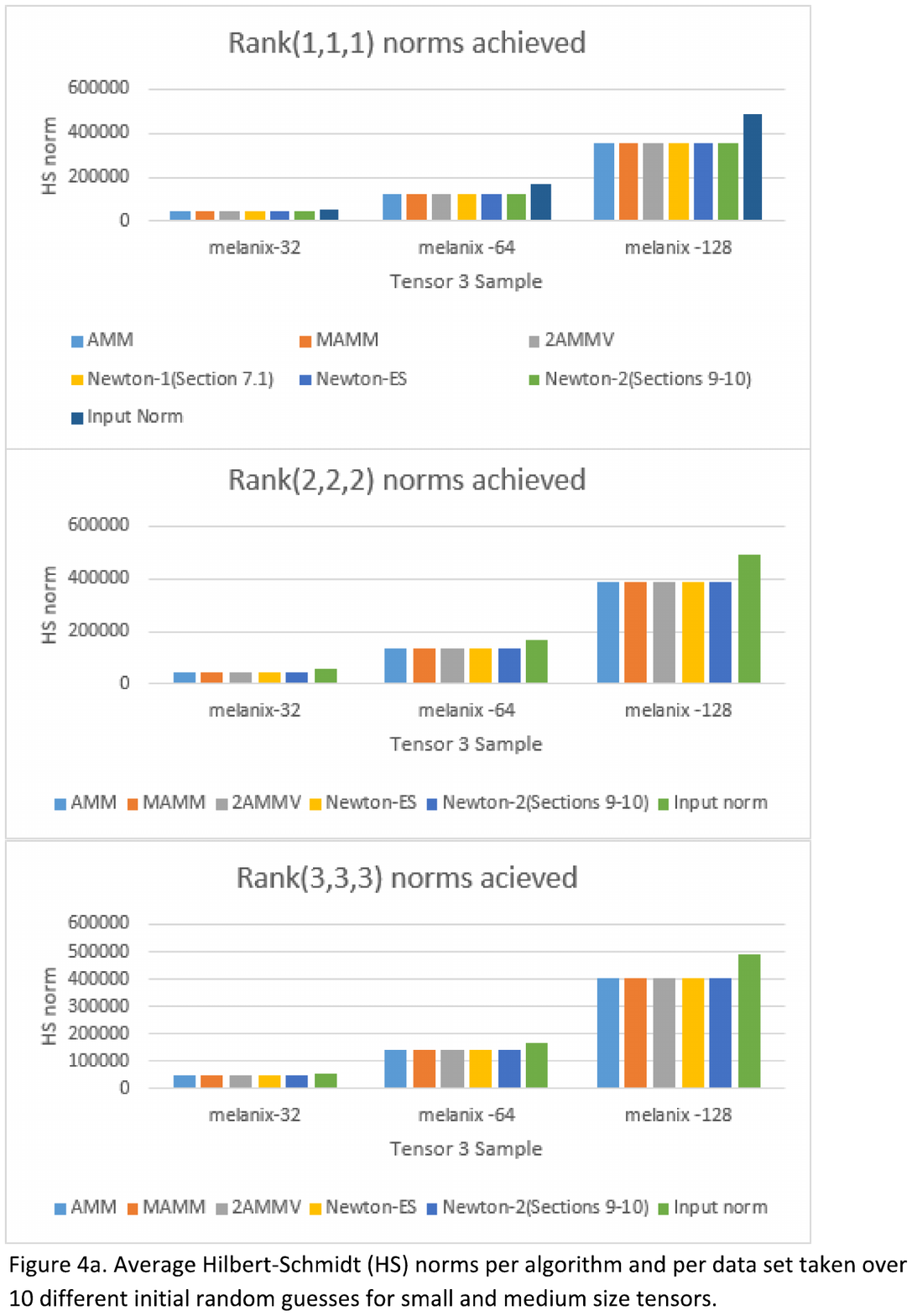}
 \includegraphics[scale = .9] {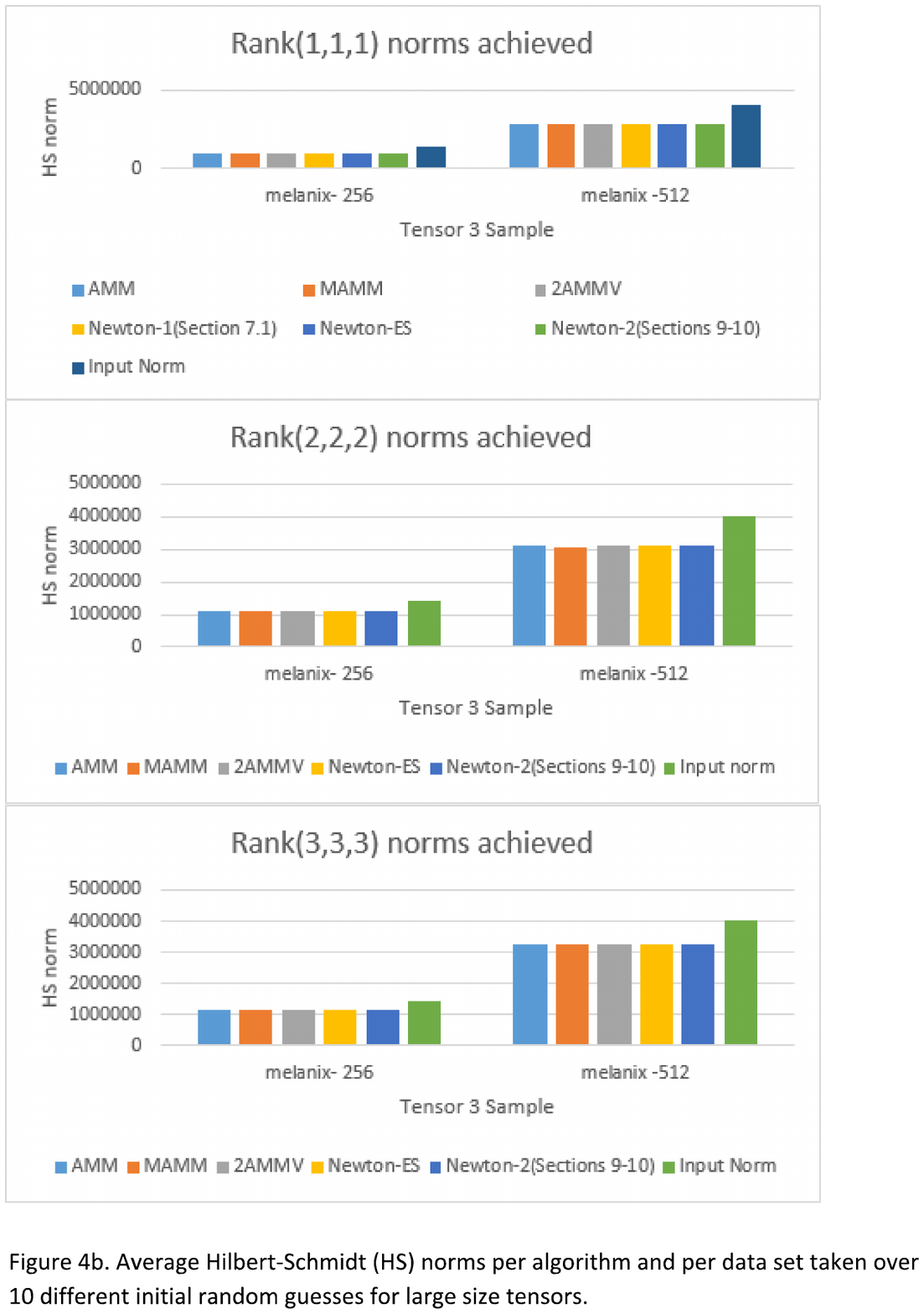}

\subsection{Best (2,2,2) and rank two approximations}\label{sub:b222apr}
Assume that $\cT$ is a $3$-tensor of rank three at least and let $\cS$ be a best $(2,2,2)$-approximation to $\cT$given by  \eqref{Sapprox}.
It is easy to show that $\cS$ has at least rank $2$.  Let $\cS'=[s_{j_1,j_2,j_3}]\in\R^{2\times 2 \times 2}$ be the core tensor corresponding 
to $\cS$.  Clearly $\rank \cS=\rank \cS'\ge 2$.
Recall that a real nonzero $2\times 2 \times 2$ tensor has rank  one, two or three \cite{BK99}.  So $\rank \cS\in\{2,3\}$.
Observe next that if $\rank\cS=\rank\cS'=2$ then $\cS$ is also a best rank two approximation of $\cT$.
Recall that a best rank two approximation of $\cT$ may not always exist.  In particular where $\rank \cT>2$ and the border rank of $\cT$ is $2$ \cite{DeSL08}.
In all our numerical simulations for best $(2,2,2)$-approximation we performed on random large tensors, the tensor $\cS'$ had rank two.  Note that the probability 
of $2\times 2\times 2$ tensors, with entries  normally distributed with mean $0$ and variance $1$, to have rank $2$ is $\frac{\pi}{4}$ \cite{Ber13}.

\section{Conclusions}\label{S:Conclusion}

We have extended the alternating maximization method (AMM) and modified alternating maximization method (MAMM) given in \cite{FMPS13}  for the computation of best rank 
one approximation to best $\br$-approximations.  We have also presented new algorithms such as $2$-alternating maximization method variant (2AMMV) and  Newton method 
for best $\br$-approximation (Newton-2).   We have provided closed form solutions for computing the $DF$ matrix in Newton-2. 
We implemented Newton-1 for best rank one approximation \cite{FMPS11} and  Newton-2. 
From the simulations, we have found out that for rank one approximation of both large and small sized tensors, Newton-1 performed the 
best.  For higher rank approximation, the best performers were 2AMMV in case of large size tensors and AMM or MAMM in case of small size tensors.\\

\emph{Acknowledgement}: We thank Daniel Kressner for his remarks. 

 \bibliographystyle{plain}
 
 \end{document}